\RequirePackage{fix-cm}
\documentclass[smallextended]{svjour3}
\smartqed 
\usepackage{graphicx}
\usepackage{amsmath,amsfonts}
\usepackage{amssymb}
\usepackage{latexsym}
\usepackage{mathrsfs}
\usepackage{colortbl}
\usepackage{amscd}
\usepackage{tikz-cd}
\usepackage{comment}
\usepackage{url}
\usepackage{color}

\allowdisplaybreaks[3]

\journalname{}

\def\diam{\mathop{\mathrm{diam}}}
\def\diag{\mathop{\mathrm{diag}}}

\def\div{\mathop{\mathrm{div}}}
\def\divh{\mathop{\mathrm{div}_h}}
\def\rot{\mathop{\mathrm{rot}}}

\def\curl{\mathop{\mathrm{curl}}}

\def\<{\mathop{\textless}}
\def\>{\mathop{\textgreater}}

\def\card{\mathop{\rm{card}}}

\spnewtheorem{thr}{Theorem}{\bf}{\it}
\spnewtheorem{defi}{Definition}{\bf}{\it}
\spnewtheorem{lem}{Lemma}{\bf}{\it}
\spnewtheorem{coro}{Corollary}{\bf}{\it}
\spnewtheorem{assume}{Assumption}{\bf}{\it}
\spnewtheorem{ex}{Example}{\it}{\it}
\spnewtheorem{Case}{Case}{\bf}{\it}
\spnewtheorem*{pf*}{Proof}{\bf}{\rm}
\spnewtheorem*{rem*}{Remark:}{\it}{\it}
\spnewtheorem*{ex*}{Example:}{\it}{\it}
\spnewtheorem{Cond}{Condition}{\bf}{\it}
\spnewtheorem{rem}{Remark}{\it}{\it}

\spnewtheorem*{lem1*}{Lemma 1}{\bf}{\rm}

\newcounter{sone}
\newcounter{stwo}
\newcounter{sthree}
\newcounter{sfour}
\newcounter{sfive}
\newcounter{ssix}
\newcounter{lone}
\newcounter{ltwo}
\newcounter{lthree}
\newcounter{lfour}
\newcounter{lfive}
\newcounter{lsix}
\makeatletter
    
    \@addtoreset{equation}{section}

\begin{document}

\title{Anisotropic modified Crouzeix--Raviart finite element method for the stationary Navier--Stokes equation 
}

\titlerunning{Anisotropic modified CR method for NS equation}  

\author{Hiroki Ishizaka 
}

\institute{Hiroki Ishizaka \at
              Team FEM, Matsuyama, Japan \\
              \email{h.ishizaka005@gmail.com}
}

\date{Received: date / Accepted: date}

\maketitle

\begin{abstract}
We studied an anisotropic modified Crouzeix--Raviart finite element method for the rotational form of a stationary incompressible Navier--Stokes equation with large irrotational body forces. We present an anisotropic $H^1$ error estimate for the velocity of the modified Crouzeix--Raviart finite element method for the Navier--Stokes equation. The modified Crouzeix--Raviart finite element scheme was obtained using a lifting operator that mapped the velocity test functions to $H(\div;\Omega)$-conforming finite element spaces. Because no shape-regularity mesh conditions are imposed, anisotropic meshes can be used for the analysis. The core idea of the proof involves using the relation between the Raviart--Thomas and Crouzeix--Raviart finite element spaces. Furthermore, we present a discrete Sobolev inequality under semi-regular mesh conditions to estimate the stability of the proposed method, and confirm the results obtained through numerical experiments.

\end{abstract}

\section{Introduction} \label{intro}
Let $\Omega \subset \mathbb{R}^d$ and $d \in \{ 2,3\}$ be a bounded polyhedral domain. The stationary Navier--Stokes problem involves determining $(u,p_0): \Omega \to \mathbb{R}^d \times \mathbb{R}$ such that
\begin{align}
\displaystyle
- \nu \varDelta u + (u \cdot \nabla ) u + \nabla p_0 = f \quad \text{in $\Omega$}, \quad \div u  = 0 \quad \text{in $\Omega$}, \quad u = 0 \quad \text{on $\partial \Omega$}, \label{intro1}
\end{align}
where $\nu$ is a nonnegative parameter and $f:\Omega \to \mathbb{R}^d$ is a given function. We use the standard notations for the Lebesgue and Sobolev spaces with associated norms \cite{ErnGue04,ErnGue21a,ErnGue21b,ErnGue21c,GirRav86,Gri11,Joh16,Soh01}. We define the curl operator \cite{GirRav86,Joh16} as  follows:
\begin{align*}
\displaystyle
\curl \varphi = \left( \frac{\partial \varphi}{\partial x_2} , - \frac{\partial \varphi}{\partial x_1} \right)^{\top}, \quad \curl u = \frac{\partial u_2}{\partial x_1} - \frac{\partial u_1}{\partial x_2}
\end{align*}
for distributions $\varphi$ of $\mathscr{D}^{\prime}(\Omega)$ and $u$ of  $ \mathscr{D}^{\prime}(\Omega)^2$ when $d=2$ and 
\begin{align*}
\displaystyle
\curl u = \left( \frac{\partial u_3}{\partial x_2} - \frac{\partial u_2}{\partial x_3} , \frac{\partial u_1}{\partial x_3} - \frac{\partial u_3}{\partial x_1} , \frac{\partial u_2}{\partial x_1} - \frac{\partial u_1}{\partial x_2} \right)^{\top}
\end{align*}
for the distribution $u$ of $\mathscr{D}^{\prime}(\Omega)^3$ when $d=3$. These operators lead to
\begin{align*}
\displaystyle
 (\curl u) \times u + \frac{1}{2} \nabla (u^{\top} u) = (u \cdot \nabla )u.
\end{align*}
This study examines the rotation form of the stationary Navier--Stokes equation \eqref{intro1}, as follows: Find $(u,p): \Omega \to \mathbb{R}^d \times \mathbb{R}$ such that
\begin{align}
\displaystyle
- \nu \varDelta u + (\curl u) \times u + \nabla p = f \quad \text{in $\Omega$}, \quad \div u  = 0 \quad \text{in $\Omega$}, \quad u = 0 \quad \text{on $\partial \Omega$}, \label{intro1b}
\end{align}
where $p$ represents the Bernoulli pressure $p := p_0 + \frac{1}{2} u^{\top} u$. 

Mixed finite element methods are commonly used in numerical analyses to solve incompressible Navier--Stokes equations. These methods require a discrete inf-sup condition to guarantee the stability and provide a unique solution. However, it is challenging to determine the Stokes elements that satisfy the discrete inf-sup condition for anisotropic mesh partitions. Furthermore, even on regular-shaped meshes, most of the Stokes elements of the inf-sup stable finite element spaces consider divergence-free discrete function spaces $V_{h,\div}$, which are non-conforming in the space of divergence-free weak functions $V_{\div}$. When the Stokes elements are used, the velocity error depends on the pressure error (\cite[Theorem 6.30, Corollary 6.33]{Joh16}); that is, 
\begin{align*}
\displaystyle
|u - u_h|_{V_h} \leq c_{1h} \inf_{v_h \in V_h} |u - v_h|_{V_h} + c_{2h} \nu^{-1} \inf_{q_h \in Q_h} \| p - q_h \|_{Q_h},
\end{align*}
where $u_h$ denotes the discrete solution of $u$. $V_h$ and $Q_h$ are the discrete velocity and pressure spaces, respectively. If $\nu$ is small and the pressure error is large, classical mixed finite element methods lead to poor velocity approximations, even when finer meshes are used. This phenomenon is known as the poor conservation of mass.

To overcome this difficulty, using the $H(\div;\Omega) := \{v \in L^2(\Omega)^d: \ \div v \in L^2(\Omega) \}$-conforming velocity lifting operator $\mathcal{L}_h$, the schemes are modified when $V_{h,\div} \not\subset V_{\div}$ by replacing the velocity test function $v_h$ in the discrete momentum equation with the lifting operator $\mathcal{L}_h(v_h)$. This idea was presented by Linke \cite{Lin14} for the Stokes equation, and has been widely applied in finite element methods \cite{Johetal17,Ledetal17,LinMatTob16,QuiPie20,YanHeZha22}. It can be said that the discretisation is pressure-robust or well balanced. In \cite{Lin14}, as the Stokes element, the Crouzeix--Raviart (CR) finite element space for velocity and the elementwise constant discontinuous space for pressure were used for the analysis. The corresponding velocity-lifting operator maps the CR velocity test functions to the lowest-order Raviart--Thomas (RT) finite element space that conforms to $H(\div;\Omega)$. The modified CR finite element method leads to a pressure-independent velocity error estimate under the shape-regular mesh condition of the form
\begin{align*}
\displaystyle
|u - u_h|_{V_h} \leq c h |u|_{H^2(\Omega)^d}.
\end{align*}

The Navier--Stokes equations have various discrete trilinear forms for the nonlinear convection term, such as convective, skew-symmetric, and (skew-symmetric) rotational forms with modified pressure (e.g., \cite{Joh16}). For continuous problems, the nonlinear convection terms do not contribute to the stability estimates. However, the discrete convective form affects the discrete stability estimates and uniqueness of the velocity. Meanwhile, because the space $H(\div;\Omega)$ is not in $H^1(\Omega)^d$, the pressure-robust schemes with skew-symmetric form are not possible for the conforming Stokes elements. In \cite{LinMer16a}, Linke and Merdon proposed a pressure-robust method for the rotation form of the Navier--Stokes equation by using the standard discrete rotation form with an additional term, which preserves skew symmetry. Subsequently, in \cite{QuiPie20}, the pressure-robust hybrid high-order (HHO) method for the skew-symmetric rotation form of the stationary Navier--Stokes equation was studied. In \cite{YanHeZha22}, the high-order pressure-robust method of the skew-symmetric rotation form of the stationary Navier--Stokes equation was proposed.

The main contributions of this study are as follows: We studied a modified CR finite element method for the skew-symmetric rotation form of the stationary incompressible Navier--Stokes equation for anisotropic mesh partitions. We present the discrete Sobolev inequality and $H^1$ error estimate of the modified CR finite element method by using a lifting operator in the lowest-order RT finite element space under relaxed mesh conditions \eqref{NewGeo}.

The remainder of this paper is organised as follows. Section 2 introduces the features of the Navier--Stokes equation in a continuous setting. Section 3 discusses the interpolation error estimates and finite element spaces. Section 4 introduces the modified CR finite element scheme and presents the main theorem of this study. Section 5 presents the numerical results. 

Throughout, we denote by $c$ a constant independent of $h$ (defined later) and of the angles and aspect ratios of simplices unless specified otherwise, and all constants $c$ are bounded {from above} if the maximum angle is bounded {from above}. These values vary across different contexts.

\section{Continuous problem}
\subsection{Weak formulation}
In finite element analysis, we derive a variational form of the continuous problem \eqref{intro1b}. We set $X := H^1(\Omega)^d$, $V := H_0^1(\Omega)^d$, and $Q := L^2_0(\Omega) := \left \{ q \in L^2(\Omega): \ \int_{\Omega} q dx = 0 \right \}$ {with a norm $\| \cdot \|_X := \| \cdot \|_{H^1(\Omega)^d}$. a seminorm $| \cdot |_V := | \cdot |_{H^1(\Omega)^d}$, and a norm $\| \cdot \|_Q := \| \cdot \|_{L^2(\Omega)}$, respectively.} The rotational form of the nonlinear term is defined as
\begin{align*}
\displaystyle
a_{\rot}(u,v,w) &:= \int_{\Omega} (\curl u) \times v \cdot w dx \quad \forall u,v,w \in X.
\end{align*}
Note that $a_{\rot}(u,v,v) = 0$ for any $u,v \in X$. $a (_\cdot,_\cdot) : X \times X \to \mathbb{R}$ and $b(_\cdot,_\cdot):X \times L^2(\Omega) \to \mathbb{R}$ denote the bilinear forms defined as
\begin{align*}
\displaystyle
a(v,w) &:= \int_{\Omega} \nabla v : \nabla w dx = \sum_{i=1}^d \int_{\Omega} \nabla v_i \cdot \nabla w_i dx, \\
b(v, q) &:= - \int_{\Omega} \div v q dx,
\end{align*}
for any $(v,w) \in X \times X$ and $(v,q) \in X \times L^2(\Omega)$. For any $f \in L^2(\Omega)^d$, the variational formulation for the Navier--Stokes equations \eqref{intro1b} is to find $(u,p) \in V \times Q$ such that
\begin{subequations} \label{ns1}
\begin{align}
\displaystyle
\nu a(u,v) + a_{\rot}(u,u,v) + b(v, p) &= \int_{\Omega} f \cdot v dx \quad \forall v \in V, \label{ns1a} \\
b(u , q) &= 0 \quad \forall q \in Q. \label{ns1b}
\end{align}
\end{subequations}
The space of the divergence-free weak function is defined as $V_{\div} := \{ v \in V: \ b(v,q) = 0 \ \forall q \in Q \}$. The reduced problem in \eqref{ns1} is as follows. Find $u \in V_{\div}$ such that
\begin{align}
\displaystyle
\nu a(u,v) + a_{\rot}(u,u,v)  &= \int_{\Omega} f \cdot v dx \quad \forall v \in V_{\div}. \label{ns2}
\end{align}

\subsection{Helmholtz projection}
We set $H_*^1(\Omega) := H^1(\Omega) \cap L^2_0(\Omega)$ and $\mathcal{H}
 = \{ v \in L^2(\Omega)^d: \ \div v = 0, \ v|_{\partial \Omega} \cdot n = 0 \}$, where $\div v = 0$ and $v|_{\partial \Omega} \cdot n = 0$ signifies that $\int_{\Omega} (v \cdot \nabla ) q dx = 0$ for any $q \in H_*^1(\Omega) $. Then, each $f \in L^2(\Omega)^d$ has a unique decomposition
\begin{align*}
\displaystyle
f = f_0 + \nabla \phi, \quad f_0 \in \mathcal{H}, \quad \phi \in H_*^1(\Omega),
\end{align*}
For example, see \cite[Lemma 74.1]{ErnGue21c} and \cite{GirRav86,Joh16,Soh01}. The $L^2$-orthogonal projection $P_{\mathcal{H}}: L^2(\Omega)^d \to \mathcal{H}$ that results from this decomposition is defined as $P_{\mathcal{H}} f := f_0$. 

Let $f \in L^2(\Omega)^d$. For any $v \in H^1(\Omega)^d$ that is divergence-free and vanishes at the boundary, it holds that
\begin{align}
\displaystyle
\int_{\Omega} f \cdot v dx &= \int_{\Omega} P_{\mathcal{H}} f \cdot v dx. \label{helm}
\end{align}
The Poincar\'e inequality yields:
\begin{align}
\displaystyle
\left|  \int_{\omega} P_{\mathcal{H}} f \cdot v dx \right| &\leq \| P_{\mathcal{H}} f \|_{L^2(\Omega)^d} \| v \|_{L^2(\Omega)^d} \leq C_P \| P_{\mathcal{H}} f \|_{L^2(\Omega)^d} | v |_{V}, \label{helm2}
\end{align}
where $C_P$ denotes the Poincar\'e constant. Furthermore, for any $u \in W^{1,4}(\Omega)^d$, $ (\curl u) \times u \in L^2(\Omega)^d$ holds. Then,
\begin{align}
\displaystyle
a_{\rot}(u,u,v) = \int_{\Omega} (\curl u) \times u \cdot v dx &=  \int_{\Omega} P_{\mathcal{H}}\left(  (\curl u) \times u \right) \cdot v dx. \label{helm1}
\end{align}

\subsection{Existence, uniqueness, and stability of a solution}
This section presents notable points regarding the continuity problem. These are the key features in the design of finite element schemes.

For the trilinear form $a_{\rot}$, from Sobolev embedding theorem $H^1(\Omega) \hookrightarrow L^4(\Omega)$, and H\"older inequality with exponents $(\frac{1}{4},\frac{1}{2},\frac{1}{4})$ we obtain
\begin{align}
\displaystyle
|a_{\rot}(u,v,w)| \leq N |u|_V |v|_V |w|_V \quad \forall u,v,w \in V, \label{triest}
\end{align}
where $N$ is a positive constant independent of $u$, $v$, and $w$. In addition, the bilinear form $a$ is continuous on $V \times V$ and coercive on $V_{\div} \times V_{\div}$. The bilinear form $b$ is continuous on $V \times Q$, and there exists a positive constant $\beta$ such that
\begin{align}
\displaystyle
\inf_{q \in Q} \sup_{v \in V} \frac{b(v , q)}{ | v |_{V} \| q \|_{Q} }\geq \beta \>0. \label{ns3}
\end{align}
For example, refer to \cite[Theorem 3.46]{Joh16}, \cite[Lemma 53.9]{ErnGue21b}, and \cite[Lemma 4.1]{GirRav86}. 
\begin{thr}
For any $f \in L^2(\Omega)^d$, there exists at least one solution of the variational formulation \eqref{ns1} for the Navier--Stokes equations.
\end{thr}

\begin{pf*}
The proof of its existence is found in \cite[Chap. \Roman{lfour}. Theorem 2.1]{GirRav86} and \cite[Chap. \Roman{ltwo}, Theorem 1.2]{Tem01}.
\qed	
\end{pf*}

We are interested in the uniqueness of the velocity solutions. This result is slightly different from the usual results. 

\begin{thr}
For $f \in L^2(\Omega)^d$, if 
\begin{align}
\displaystyle
\frac{N C_P}{\nu^2} \| P_{\mathcal{H}} f \|_{L^2(\Omega)^d} \< 1 \label{ns6}
\end{align}
holds true, problem \eqref{ns2} has a unique solution $u \in V_{\div}$.
\end{thr}

\begin{pf*}
We provide a proof of the uniqueness as follows. Let $u \in V_{\div}$ be the solution to problem \eqref{ns2}. By setting $v := u$ in \eqref{ns2} and using \eqref{helm2}, we obtain
\begin{align*}
\displaystyle
\nu |u|^2_V = \int_{\omega} P_{\mathcal{H}} f \cdot u dx \leq C_P \| P_{\mathcal{H}} f \|_{L^2(\Omega)^d} | u |_{V},
\end{align*}
which leads to
\begin{align}
\displaystyle
|u|_V \leq \frac{C_P}{\nu}\| P_{\mathcal{H}} f \|_{L^2(\Omega)^d}, \label{ns7}
\end{align}
where $C_P$ denotes the Poincar\'e constant. Let $u_1$ and $u_2$ be solutions to problem \eqref{ns2}. We set $u := u_1 - u_2$ and subtract the problem \eqref{ns2} corresponding to $u_1$ and $u_2$. Then,
\begin{align}
\displaystyle
\nu a(u,v) + a_{\rot}(u_2,u,v) + a_{\rot}(u,u_1,v) = 0 \quad \forall v \in V_{\div}. \label{ns8}
\end{align}
Setting $v := u$ in \eqref{ns8} with \eqref{triest} and \eqref{ns7}, it holds that
\begin{align*}
\displaystyle
\nu |u|^2_V 
&= - a_{\rot}(u,u_1,u) 
\leq N |u|^2_V |u_1|_V \leq \frac{N C_P}{\nu} \| P_{\mathcal{H}} f \|_{L^2(\Omega)^d}  |u|^2_V,
\end{align*}
which leads to
\begin{align*}
\displaystyle
\left(1 -   \frac{N C_P}{\nu^2} \| P_{\mathcal{H}} f \|_{L^2(\Omega)^d}  \right) |u|^2_V \leq 0.
\end{align*}
If \eqref{ns6} holds true, this inequality implies that $u_1 = u_2$.
\qed
\end{pf*}

\begin{rem}
The stability estimate \eqref{ns7} is sharper than the classical estimate $|u|_V \leq \frac{C_p}{\nu} \| f \|_{L^2(\Omega)^s}$. The stability estimate for the pressure follows from the inf-sup condition \eqref{ns3}. Changing the right-hand side of \eqref{ns2} to $f \to f + \nabla \phi$ for $\phi \in H_*^1(\Omega)$ yields 
\begin{align*}
\displaystyle
\nu |u|_V \leq C_P \| P_{\mathcal{H}} (f + \nabla \phi) \|_{L^2(\Omega)^d} = C_P \| P_{\mathcal{H}} f \|_{L^2(\Omega)^d},
\end{align*}
as $P_{\mathcal{H}}(\nabla \phi) = 0$. Therefore, the velocity is unaffected by the irrotational part of the force.
\end{rem}

\section{Modified CR finite element method}

\subsection{Notation}
For $k \in \mathbb{N}_0 := \mathbb{N} \cup \{ 0 \}$, $\mathbb{P}^k(T)$ is spanned by the restriction to $T$ by the polynomials in $\mathbb{P}^k$, where  $\mathbb{P}^k$ denotes the space of polynomials with a maximum of $k$ degrees.

\textbf{Meshes, mesh faces, and jumps.} Let $\mathbb{T}_h = \{ T \}$ be a simplicial mesh of $\overline{\Omega}$ composed of closed $d$-simplices:
\begin{align*}
\displaystyle
\overline{\Omega} = \bigcup_{T \in \mathbb{T}_h} T,
\end{align*}
where $h := \max_{T \in \mathbb{T}_h} h_{T}$ and $ h_{T} := \diam(T)$. For simplicity, we assume that $\mathbb{T}_h$ is conformal; that is, $\mathbb{T}_h$ is a simplicial mesh of $\overline{\Omega}$ without hanging nodes.

Let $\mathcal{F}_h^i$ be the set of interior faces, and $\mathcal{F}_h^{\partial}$ be the set of faces on boundary $\partial \Omega$. We set $\mathcal{F}_h := \mathcal{F}_h^i \cup \mathcal{F}_h^{\partial}$. For any $F \in \mathcal{F}_h$, we define the unit normal $n_F$ to $F$ as follows: (\roman{sone}) If  $F \in \mathcal{F}_h^i$ with $F = T_1 \cap T_2$, $T_1,T_2 \in \mathbb{T}_h$, let $n_1$ and $n_2$ be the outwards unit normals of $T_1$ and $T_2$, respectively. Then, $n_F$ is either of $\{ n_1 , n_2\}$; (\roman{stwo}) If $F \in \mathcal{F}_h^{\partial}$, $n_F$ is the unit outwards normal $n$ to $\partial \Omega$. {We define the broken (piecewise) Sobolev space as
\begin{align*}
\displaystyle
H^1(\mathbb{T}_h) &:= \left\{ \varphi \in L^2(\Omega); \ \varphi|_{T} \in H^1(T) \ \forall T \in \mathbb{T}_h  \right\}
\end{align*}
with the seminorm
\begin{align*}
\displaystyle
| \varphi |_{H^1(\mathbb{T}_h)} &:= \left( \sum_{T \in \mathbb{T}_h} | \varphi |^2_{H^1(T)} \right)^{\frac{1}{2}} \quad \varphi \in H^1(\mathbb{T}_h).
\end{align*}
}
Let $\varphi \in H^1(\mathbb{T}_h)$. Suppose that $F \in \mathcal{F}_h^i$ with $F = T_1 \cap T_2$, $T_1,T_2 \in \mathbb{T}_h$. We set $\varphi_1 := \varphi{|_{T_1}}$ and $\varphi_2 := \varphi{|_{T_2}}$. We set two nonnegative real numbers $\omega_{T_1,F}$ and $\omega_{T_2,F}$ such that
\begin{align*}
\displaystyle
\omega_{T_1,F} + \omega_{T_2,F} = 1.
\end{align*}
The jump and skew-weighted averages of $\varphi$ across $F$ are then defined as
\begin{align*}
\displaystyle
[\![\varphi]\!] := [\! [ \varphi ]\!]_F := \varphi_1 - \varphi_2, \quad  \{\! \{ \varphi\} \! \}_{\overline{\omega}} :=  \{\! \{ \varphi\} \! \}_{\overline{\omega},F} := \omega_{T_2,F} \varphi_1 + \omega_{T_1,F} \varphi_2.
\end{align*}
For a boundary face $F \in \mathcal{F}_h^{\partial}$ with $F = \partial T \cap \partial \Omega$, $[\![\varphi ]\!]_F := \varphi|_{T}$ and $\{\! \{ \varphi \} \!\}_{\overline{\omega}} := \varphi |_{T}$. For any $v \in H^{1}(\mathbb{T}_h)^d$, the notation
\begin{align*}
\displaystyle
&[\![v \cdot n]\!] := [\![ v \cdot n ]\!]_F := v_1 \cdot n_F - v_2 \cdot n_F,  \  \{\! \{ v\} \! \}_{\omega} :=  \{\! \{ v \} \! \}_{\omega,F} := \omega_{T_1,F} v_1 + \omega_{T_2,F} v_2
\end{align*}
for the jump in the normal component and weighted average of $v$. For any $v \in H^{1}(\mathbb{T}_h)^d$ and $\varphi \in H^{1}(\mathbb{T}_h)$, we have that
\begin{align*}
\displaystyle
[\![ (v \varphi) \cdot n ]\!]_F
&=  \{\! \{ v \} \! \}_{\omega,F} \cdot n_F [\! [ \varphi ]\!]_F + [\![ v \cdot n ]\!]_F \{\! \{ \varphi\} \! \}_{\overline{\omega},F}.
\end{align*}

The broken gradient operator is defined as follows: For $\varphi \in H^{1}(\mathbb{T}_h)$, the broken gradient $\nabla_h: H^{1}(\mathbb{T}_h) \to L^2(\Omega)^{d}$ is defined as
\begin{align*}
\displaystyle
(\nabla_h \varphi)|_{T} &:= \nabla (\varphi|_{T}) \quad \forall T \in \mathbb{T}_h,
\end{align*}
Furthermore, we define a broken $H(\div;T)$-space as
\begin{align*}
\displaystyle
H(\div;\mathbb{T}_h) := \left \{ v \in L^2(\Omega)^d; \ v |_{T} \in H(\div;T) \ \forall T \in \mathbb{T}_h  \right\}.
\end{align*}
Thus, the broken divergence operator $\divh : H(\div;\mathbb{T}_h) \to L^2(\Omega)$, such that for all $v \in H(\div;\mathbb{T}_h)$,
\begin{align*}
\displaystyle
(\divh v)|_{T} := \div (v |_{T}) \quad \forall T \in \mathbb{T}_h.
\end{align*}

\subsection{Finite element spaces} \label{FEspaces0}
For $m \in \mathbb{N}_0$, we define a discontinuous finite element space as
\begin{align*}
\displaystyle
P_{dc,h}^m &:= \left\{ p_h \in L^{\infty}(\Omega); \ p_h|_{T} \in \mathbb{P}^{m}({T}) \quad \forall T \in \mathbb{T}_h \right\}.
\end{align*}
The CR finite element space is defined as
\begin{align*}
\displaystyle
V_{h0}^{CR} &:=  \biggl \{ \varphi_h \in P_{dc,h}^1: \  \int_F [\![ \varphi_h ]\!] ds = 0 \ \forall F \in \mathcal{F}_h \biggr \}.
\end{align*}
Thus, we define the Stokes pair $(V_h,Q_h)$ as
\begin{align*}
\displaystyle
V_h := (V_{h0}^{CR})^d, \quad Q_h := P_{dc,h}^0 \cap L^{2}_0(\Omega)
\end{align*}
with norms
\begin{align*}
\displaystyle
|v_h|_{V_{h}} :=  |v_h|_{H^{1}(\mathbb{T}_h)^d} = \left ( \sum_{i=1}^d |v_{h,i}|^2_{H^{1}(\mathbb{T}_h)} \right)^{1/2}, \quad \| q_h \|_{Q_h} := \| q_h \|_{L^2(\Omega)}
\end{align*}
for any $v_h = (v_{h,1},\ldots,v_{h,d})^{\top}  \in V_{h}$ and $q_h \in Q_h$. 

\subsection{Approximation of the trilinear form $a_{\rot}$}
The following lemma generalises \cite[Lemma 6.7]{Joh16} (see also \cite{QuiPie20,YanHeZha22}).

\begin{lem} \label{iden=lem5}
 It holds that for any $u,v,w \in H^1(\mathbb{T}_h)^d$
 \begin{align}
\displaystyle
\sum_{T \in \mathbb{T}_h} \int_{T} (\curl u) \times v \cdot w dx  &=  \sum_{T \in \mathbb{T}_h} \int_{T} \{ (v \cdot \nabla) u \cdot w - (w \cdot \nabla )u \cdot v \} dx.  \label{rot=uvw}
\end{align}
\end{lem}

\begin{pf*}
Let $u = (u_1,\ldots,u_d)^{\top} \in H^1(\mathbb{T}_h)^d$, $v = (v_1,\ldots,v_d)^{\top} \in H^1(\mathbb{T}_h)^d$, and $w = (w_1,\ldots,w_d)^{\top} \in H^1(\mathbb{T}_h)^d$. The following identity holds:
\begin{align*}
\displaystyle
(\curl u) \times v
= (v \cdot \nabla) u - \nabla (u \cdot v) + (\nabla v)^{\top} u,
\end{align*}
where $\displaystyle \nabla v := \left(\frac{\partial v_i}{\partial x_j} \right)_{i,j}$. For any $T \in \mathbb{T}_h$, we obtain
\begin{align}
\displaystyle
\int_T (\curl u) \times v \cdot w dx = \int_T \left\{ (v \cdot \nabla) u - \nabla (u \cdot v) + (\nabla v)^{\top} u \right\} \cdot w dx. \label{iden4=2}
\end{align}
From Green's formula, we obtain
\begin{align}
\displaystyle
- \int_T  \nabla (u \cdot v) \cdot w dx &= - \int_{\partial T} (u \cdot v) (w \cdot n_T) ds + \int_{T} (u \cdot v) \div w dx, \label{iden4=3}\\
\int_T (\nabla v)^{\top} u \cdot w dx
&= \int_T \nabla v : u w^{\top} dx \notag\\
&\hspace{-2cm} = \int_{\partial T} (u \cdot v) (w \cdot n_T) ds - \int_{T} v \cdot \div(u w^{\top}) dx \notag \\
&\hspace{-2cm} = \int_{\partial T} (u \cdot v) (w \cdot n_T) ds - \int_T (u \cdot v) \div w dx - \int_T (w \cdot \nabla u) \cdot v dx. \label{iden4=4}
\end{align}
From \eqref{iden4=2}, \eqref{iden4=3}, and \eqref{iden4=4}, it holds that for any $T \in \mathbb{T}_h$,
\begin{align*}
\displaystyle
\int_T (\curl u) \times v \cdot w dx =  \int_{T} \{ (v \cdot \nabla) u \cdot w - (w \cdot \nabla )u \cdot v \} dx.
\end{align*}
By summing over $T$, the target identity \eqref{rot=uvw} is obtained.
\qed
\end{pf*}

%\begin{rem}
%The equality \eqref{rot=uvw} is valid for any $u,v \in H^1(\mathbb{T}_h)^d$ and $w \in H(\div;\mathbb{T}_h)$.
%\end{rem}

For the classic CR finite element methods, the discretisation of the nonlinear form is as follows:
\begin{align*}
\displaystyle
\sum_{T \in \mathbb{T}_h} \int_{T}  (\curl u_h) \times v_h \cdot w_h dx \quad \forall u_h,v_h,w_h \in  V_h.
\end{align*}
As a lifting operator, we consider the RT interpolation operator $\mathcal{I}_{h}^{RT^0}: V_h  \cup W^{1,1}(\Omega)^d \to V^{RT^0}_{h}$ defined in Section \ref{RTsp}. The associated  discrete convective trilinear form  to the trilinear form \eqref{helm1} is defined as follows:
\begin{align*}
\displaystyle
\sum_{T \in \mathbb{T}_h} \int_{T}  (\curl u_h) \times u_h \cdot \mathcal{I}_{h}^{RT^0} v_h dx \quad \forall u_h,v_h \in  V_h.
\end{align*}
Applying Lemma \ref{iden=lem5} yields
\begin{align}
\displaystyle
\sum_{T \in \mathbb{T}_h} \int_T \{ ( u_h \cdot \nabla) u_h \cdot \mathcal{I}_h^{RT^0} v_h - (\mathcal{I}_h^{RT^0} v_h \cdot \nabla )u_h \cdot  u_h \} dx. \label{conv=term}
\end{align}
However, this form does not vanish by substituting $u_h$ for $v_h$. Therefore, we define a discrete trilinear form $a_{\rot,h}: (V+V_h) \times (V+V_h) \times (V+V_h) \to \mathbb{R}$ as
\begin{align*}
\displaystyle
&a_{\rot,h}(u_h, v_h, w_h ) \\
&\quad := \sum_{T \in \mathbb{T}_h} \int_T \left \{ (\mathcal{I}_h^{RT^0} v_h \cdot \nabla) u_h \cdot \mathcal{I}_h^{RT^0} w_h - (\mathcal{I}_h^{RT^0} w_h \cdot \nabla )u_h \cdot \mathcal{I}_h^{RT^0} v_h \right \} dx
\end{align*}
for any $u_h,v_h, w_h \in V + V_h$. The trilinear form $a_{\rot,h}$ has two properties. For any $u_h,v_h,w_h \in V_h$,
\begin{align*}
\displaystyle
&a_{\rot,h}(u_h, v_h, v_h ) = 0. \\
&a_{\rot,h}(u_h, v_h, w_h ) = - a_{\rot,h}(u_h, w_h, v_h ).
\end{align*}

\subsection{Modified CR finite element scheme}
We consider the CR finite element method for the Navier--Stokes equation \eqref{ns1} as follows: Find $(u_h,p_h) \in V_{h} \times Q_h$ such that
\begin{subequations} \label{cr=1}
\begin{align}
\displaystyle
\nu a_h(u_h,v_h) + a_{\rot,h}(u_h,  u_h,  v_h ) + b_h(v_h , p_h)
&= \int_{\Omega} f \cdot \mathcal{I}_h^{RT^0} v_h dx \quad \forall v_h \in V_{h}, \label{cr=1a} \\
b_h(u_h , q_h) &= 0 \quad \forall q_h \in Q_h, \label{cr=1b}
\end{align}
\end{subequations}
where the discrete bilinear forms $a_{h}: (V+ V_{h}) \times (V + V_{h}) \to \mathbb{R}$ and $b_h: (V +V_{h}) \times Q_h \to \mathbb{R}$ are the discrete counterparts of the bilinear forms $a$ and $b$ defined as
\begin{align*}
\displaystyle
a_{h}(u_h,v_h) &:= \sum_{i=1}^d \sum_{T \in \mathbb{T}_h} \int_{T} \nabla u_{h,i} \cdot \nabla v_{h,i} dx, \\
b_h(v_h , q_h) &:= - \sum_{T \in \mathbb{T}_h} \int_{T} \div v_h q_h dx.
\end{align*}

We define a discrete divergence-free weak subspace as follows:
\begin{align*}
\displaystyle
V_{h,\div} := \{ v_h \in V_h: \ b_h(v_h,q_h) = 0 \ \forall q_h \in Q_h \}.
\end{align*}
Because $V_{h,\div} \not\subset V_{\div}$, the space $V_{h,\div}$ is nonconforming in the space $V_{\div}$. Then, the reduced problem of \eqref{cr=1} is as follows. Find $u_h \in V_{h,\div}$ such that
\begin{align}
\displaystyle
\nu a_{h}(u_h,v_h) + a_{\rot,h}( u_h, u_h, v_h ) = \int_{\Omega} f \cdot \mathcal{I}_{h}^{RT^0} v_h dx \quad \forall v_h \in V_{h,\div}. \label{cr=2}
\end{align}

The bilinear form $a_h$ is continuous and coercive on $V_h \times V_h$, and the bilinear form $b_h$ is continuous on $V_h \times Q_h$ and satisfies the inf-sup condition; there exists a positive constant $\beta_0$ such that
\begin{align}
\displaystyle
\inf_{q_h \in Q_h} \sup_{v_h \in V_h} \frac{|b_h(v_h,q_h)|}{|v_h|_{V_h} \| q_h \|_{Q_h} } \geq \beta_0 \> 0, \label{cr=4}
\end{align}
because the CR interpolation operator acts as a nonconforming Fortin operator.

\section{Anisotropic interpolation error estimates}
Our strategy for anisotropic interpolation errors in anisotropic meshes was proposed by \cite{Ish21,Ish23a,IshKobTsu23a}.

\subsection{Reference elements} \label{reference}
We now define the reference elements $\widehat{T} \subset \mathbb{R}^d$.

\subsubsection*{Two-dimensional case} \label{reference2d}
Let $\widehat{T} \subset \mathbb{R}^2$ be a reference triangle with vertices $\hat{p}_1 := (0,0)^{\top}$, $\hat{p}_2 := (1,0)^{\top}$, and $\hat{p}_3 := (0,1)^{\top}$. 

\subsubsection*{Three-dimensional case} \label{reference3d}
In the three-dimensional case, we consider the following two cases: (\roman{sone}) and (\roman{stwo}) (see Condition \ref{cond2}).

Let $\widehat{T}_1$ and $\widehat{T}_2$ be reference tetrahedra with the following vertices:
\begin{description}
   \item[(\roman{sone})] $\widehat{T}_1$ has the vertices $\hat{p}_1 := (0,0,0)^{\top}$, $\hat{p}_2 := (1,0,0)^{\top}$, $\hat{p}_3 := (0,1,0)^{\top}$, and $\hat{p}_4 := (0,0,1)^{\top}$;
 \item[(\roman{stwo})] $\widehat{T}_2$ has the vertices $\hat{p}_1 := (0,0,0)^{\top}$, $\hat{p}_2 := (1,0,0)^{\top}$, $\hat{p}_3 := (1,1,0)^{\top}$, and $\hat{p}_4 := (0,0,1)^{\top}$.
\end{description}
Therefore, we set $\widehat{T} \in \{ \widehat{T}_1 , \widehat{T}_2 \}$. The case (\roman{sone}) is called the \textit{regular vertex property}; see \cite{AcoDur99}.

\subsection{Two-step affine mapping} \label{two=step}
To an affine simplex $T \subset \mathbb{R}^d$, we construct two affine mappings $\Phi_{\widetilde{T}}: \widehat{T} \to \widetilde{T}$ and $\Phi_{T}: \widetilde{T} \to T$. First, we define the affine mapping $\Phi_{\widetilde{T}}: \widehat{T} \to \widetilde{T}$ as
\begin{align}
\displaystyle
\Phi_{\widetilde{T}}: \widehat{T} \ni \hat{x} \mapsto \tilde{x} := \Phi_{\widetilde{T}}(\hat{x}) := {A}_{\widetilde{T}} \hat{x} \in  \widetilde{T}, \label{aff=1}
\end{align}
where ${A}_{\widetilde{T}} \in \mathbb{R}^{d \times d}$ is an invertible matrix. We then define the affine mapping $\Phi_{T}: \widetilde{T} \to T$ as follows:
\begin{align}
\displaystyle
\Phi_{T}: \widetilde{T} \ni \tilde{x} \mapsto x := \Phi_{T}(\tilde{x}) := {A}_{T} \tilde{x} + b_{T} \in T, \label{aff=2}
\end{align}
where $b_{T} \in \mathbb{R}^d$ is a vector and ${A}_{T} \in O(d)$ denotes the rotation and mirror-imaging matrix. We define the affine mapping $\Phi: \widehat{T} \to T$ as
\begin{align*}
\displaystyle
\Phi := {\Phi}_{T} \circ {\Phi}_{\widetilde{T}}: \widehat{T} \ni \hat{x} \mapsto x := \Phi (\hat{x}) =  ({\Phi}_{T} \circ {\Phi}_{\widetilde{T}})(\hat{x}) = {A} \hat{x} + b_{T} \in T, 
\end{align*}
where ${A} := {A}_{T} {A}_{\widetilde{T}} \in \mathbb{R}^{d \times d}$.

\subsubsection*{Construct mapping $\Phi_{\widetilde{T}}: \widehat{T} \to \widetilde{T}$} \label{sec221} 
We consider the affine mapping \eqref{aff=1}. We define the matrix $ {A}_{\widetilde{T}} \in \mathbb{R}^{d \times d}$ as follows. We first define the diagonal matrix as
\begin{align}
\displaystyle
\widehat{A} :=  \diag (h_1,\ldots,h_d), \quad h_i \in \mathbb{R}_+ \quad \forall i,\label{aff=3}
\end{align}
where $\mathbb{R}_+$ denotes the set of positive real numbers.

For $d=2$, we define the regular matrix $\widetilde{A} \in \mathbb{R}^{2 \times 2}$ as
\begin{align}
\displaystyle
\widetilde{A} :=
\begin{pmatrix}
1 & s \\
0 & t \\
\end{pmatrix}, \label{aff=4}
\end{align}
with the parameters
\begin{align*}
\displaystyle
s^2 + t^2 = 1, \quad t \> 0.
\end{align*}
For the reference element $\widehat{T}$, let $\mathfrak{T}^{(2)}$ be a family of triangles.
\begin{align*}
\displaystyle
\widetilde{T} &= \Phi_{\widetilde{T}}(\widehat{T}) = {A}_{\widetilde{T}} (\widehat{T}), \quad {A}_{\widetilde{T}} := \widetilde {A} \widehat{A}
\end{align*}
with the vertices $\tilde{p}_1 := (0,0)^{\top}$, $\tilde{p}_2 := (h_1,0)^{\top}$ and $\tilde{p}_3 :=(h_2 s , h_2 t)^{\top}$. Then, $h_1 = |\tilde{p}_1 - \tilde{p}_2| \> 0$ and $h_2 = |\tilde{p}_1 - \tilde{p}_3| \> 0$. 

For $d=3$, we define the regular matrices $\widetilde{A}_1, \widetilde{A}_2 \in \mathbb{R}^{3 \times 3}$ as follows:
\begin{align}
\displaystyle
\widetilde{A}_1 :=
\begin{pmatrix}
1 & s_1 & s_{21} \\
0 & t_1  & s_{22}\\
0 & 0  & t_2\\
\end{pmatrix}, \
\widetilde{A}_2 :=
\begin{pmatrix}
1 & - s_1 & s_{21} \\
0 & t_1  & s_{22}\\
0 & 0  & t_2\\
\end{pmatrix} \label{aff=5}
\end{align}
with the parameters
\begin{align*}
\displaystyle
\begin{cases}
s_1^2 + t_1^2 = 1, \ s_1 \> 0, \ t_1 \> 0, \ h_2 s_1 \leq h_1 / 2, \\
s_{21}^2 + s_{22}^2 + t_2^2 = 1, \ t_2 \> 0, \ h_3 s_{21} \leq h_1 / 2.
\end{cases}
\end{align*}
Therefore, we set $\widetilde{A} \in \{ \widetilde{A}_1 , \widetilde{A}_2 \}$. For the reference elements $\widehat{T}_i$, $i=1,2$, let $\mathfrak{T}_i^{(3)}$, $i=1,2$, be a family of tetrahedra.
\begin{align*}
\displaystyle
\widetilde{T}_i &= \Phi_{\widetilde{T}_i} (\widehat{T}_i) =  {A}_{\widetilde{T}_i} (\widehat{T}_i), \quad {A}_{\widetilde{T}_i} := \widetilde {A}_i \widehat{A}, \quad i=1,2,
\end{align*}
with the vertices
\begin{align*}
\displaystyle
&\tilde{p}_1 := (0,0,0)^{\top}, \ \tilde{p}_2 := (h_1,0,0)^{\top}, \ \tilde{p}_4 := (h_3 s_{21}, h_3 s_{22}, h_3 t_2)^{\top}, \\
&\begin{cases}
\tilde{p}_3 := (h_2 s_1 , h_2 t_1 , 0)^{\top} \quad \text{for case (\roman{sone})}, \\
\tilde{p}_3 := (h_1 - h_2 s_1, h_2 t_1,0)^{\top} \quad \text{for case (\roman{stwo})}.
\end{cases}
\end{align*}
Subsequently, $h_1 = |\tilde{p}_1 - \tilde{p}_2| \> 0$, $h_3 = |\tilde{p}_1 - \tilde{p}_4| \> 0$, and
\begin{align*}
\displaystyle
h_2 =
\begin{cases}
|\tilde{p}_1 - \tilde{p}_3| \> 0  \quad \text{for case (\roman{sone})}, \\
|\tilde{p}_2 - \tilde{p}_3| \> 0  \quad \text{for case (\roman{stwo})}.
\end{cases}
\end{align*}

\subsubsection*{Construct mapping $\Phi_{T}: \widetilde{T} \to T$}  \label{sec322}
We determine the affine mapping \eqref{aff=2} as follows. Let ${T} \in \mathbb{T}_h$ have vertices ${p}_i$ ($i=1,\ldots,d+1$). Let $b_{T} \in \mathbb{R}^d$ be the vector and ${A}_{T} \in O(d)$ be the rotation and mirror imaging matrix such that
\begin{align*}
\displaystyle
p_{i} = \Phi_T (\tilde{p}_i) = {A}_{T} \tilde{p}_i + b_T, \quad i \in \{1, \ldots,d+1 \},
\end{align*}
where vertices $p_{i}$ ($i=1,\ldots,d+1$) satisfy the following conditions:

\begin{Cond}[Case in which $d=2$] \label{cond1}
Let ${T} \in \mathbb{T}_h$ have vertices ${p}_i$ ($i=1,\ldots,3$). We assume that $\overline{{p}_2 {p}_3}$ is the longest edge of ${T}$, that is, $ h_{{T}} := |{p}_2 - {p}_ 3|$. We set $h_1 = |{p}_1 - {p}_2|$ and $h_2 = |{p}_1 - {p}_3|$. We then assume that $h_2 \leq h_1$. {Because $\frac{1}{2} h_T < h_1 \leq h_T$, ${h_1 \approx h_T}$.} 
\end{Cond}

\begin{Cond}[Case in which $d=3$] \label{cond2}
Let ${T} \in \mathbb{T}_h$ have vertices ${p}_i$ ($i=1,\ldots,4$). Let ${L}_i$ ($1 \leq i \leq 6$) be the edges of ${T}$. We denote by ${L}_{\min}$  the edge of ${T}$ with the minimum length; that is, $|{L}_{\min}| = \min_{1 \leq i \leq 6} |{L}_i|$. We set $h_2 := |{L}_{\min}|$ and assume that 
\begin{align*}
\displaystyle
&\text{the endpoints of ${L}_{\min}$ are either $\{ {p}_1 , {p}_3\}$ or $\{ {p}_2 , {p}_3\}$}.
\end{align*}
Among the four edges sharing an endpoint with ${L}_{\min}$, we consider the longest edge ${L}^{({\min})}_{\max}$. Let ${p}_1$ and ${p}_2$ be the endpoints of edge ${L}^{({\min})}_{\max}$. Thus, we have
\begin{align*}
\displaystyle
h_1 = |{L}^{(\min)}_{\max}| = |{p}_1 - {p}_2|.
\end{align*}
We consider cutting $\mathbb{R}^3$ with a plane that contains the midpoint of the edge ${L}^{(\min)}_{\max}$ and is perpendicular to the vector ${p}_1 - {p}_2$. Thus, there are two cases. 
\begin{description}
  \item[(Type \roman{sone})] ${p}_3$ and ${p}_4$  belong to the same half-space;
  \item[(Type \roman{stwo})] ${p}_3$ and ${p}_4$  belong to different half-spaces.
\end{description}
In each case, we set
\begin{description}
  \item[(Type \roman{sone})] ${p}_1$ and ${p}_3$ as the endpoints of ${L}_{\min}$, that is, $h_2 =  |{p}_1 - {p}_3| $;
  \item[(Type \roman{stwo})] ${p}_2$ and ${p}_3$ as the endpoints of ${L}_{\min}$, that is, $h_2 =  |{p}_2 - {p}_3| $.
\end{description}
Finally, we set $h_3 = |{p}_1 - {p}_4|$. We implicitly assume that ${p}_1$ and ${p}_4$ belong to the same half space. Additionally, note that ${h_1 \approx h_T}$.
\end{Cond}

\subsection{Piola transformations}
The Piola transformation $\Psi : L^1(\widehat{T})^d \to L^1({T})^d$ is defined as follows:
\begin{align*}
\displaystyle
\Psi :  L^1(\widehat{T})^d  &\to  L^1({T})^d \\
\hat{v} &\mapsto v(x) :=  \Psi(\hat{v})(x) = \frac{1}{\det(A)} A \hat{v}(\hat{x}).
\end{align*}

\subsection{Additional notations and assumptions} \label{addinot}
We define the vectors ${r}_n \in \mathbb{R}^d$, $n=1,\ldots,d$ as follows: If $d=2$,
\begin{align*}
\displaystyle
{r}_1 := \frac{p_2 - p_1}{|p_2 - p_1|}, \quad {r}_2 := \frac{p_3 - p_1}{|p_3 - p_1|},
\end{align*}
and if $d=3$,
\begin{align*}
\displaystyle
&{r}_1 := \frac{p_2 - p_1}{|p_2 - p_1|}, \quad {r}_3 := \frac{p_4 - p_1}{|p_4 - p_1|}, \quad
\begin{cases}
\displaystyle
{r}_2 := \frac{p_3 - p_1}{|p_3 - p_1|}, \quad \text{for (Type \roman{sone})}, \\
\displaystyle
{r}_2 := \frac{p_3 - p_2}{|p_3 - p_2|} \quad \text{for (Type \roman{stwo})}.
\end{cases}
\end{align*}
For a sufficiently smooth function $\varphi$ and a vector function $v := (v_{1},\ldots,v_{d})^{\top}$, we define the directional derivative for $i \in \{ 1, \ldots,d \}$ as
\begin{align*}
\displaystyle
\frac{\partial \varphi}{\partial {r_i}} &:= ( {r}_i \cdot  \nabla_{x} ) \varphi = \sum_{i_0=1}^d ({r}_i)_{i_0} \frac{\partial \varphi}{\partial x_{i_0}^{}}, \\
\frac{\partial v}{\partial r_i} &:= \left(\frac{\partial v_{1}}{\partial r_i}, \ldots, \frac{\partial v_{d}}{\partial r_i} \right)^{\top}
= ( ({r}_i  \cdot \nabla_{x}) v_{1}, \ldots, ({r}_i  \cdot \nabla_{x} ) v_{d} )^{\top}.
\end{align*}
For a multiindex $\beta = (\beta_1,\ldots,\beta_d) \in \mathbb{N}_0^d$, we use the notation
\begin{align*}
\displaystyle
\partial^{\beta} \varphi := \frac{\partial^{|\beta|} \varphi}{\partial x_1^{\beta_1} \ldots \partial x_d^{\beta_d}}, \quad \partial^{\beta}_{r} \varphi := \frac{\partial^{|\beta|} \varphi}{\partial r_1^{\beta_1} \ldots \partial r_d^{\beta_d}}, \quad h^{\beta} :=  h_{1}^{\beta_1} \cdots h_{d}^{\beta_d}.
\end{align*}
It should be noted that $\partial^{\beta} \varphi \neq  \partial^{\beta}_{r} \varphi$.

We propose a geometric parameter $H_{T}$ in a prior paper \cite{IshKobTsu23a}.
 \begin{defi} \label{defi1}
 The parameter $H_{{T}}$ is defined as
\begin{align*}
\displaystyle
H_{{T}} := \frac{\prod_{i=1}^d h_i}{|{T}|_d} h_{{T}},
\end{align*}
where {$|\cdot|_d$ denotes the $d$-dimensional Hausdorff measure.}

\end{defi}
Here, we introduce the geometric condition proposed by \cite{IshKobTsu21a,IshKobTsu23a}, which is equivalent to the maximum angle condition \cite{IshKobSuzTsu21d,IshKobTsu23a}.

\begin{assume} \label{neogeo=assume}
A family of meshes $\{ \mathbb{T}_h\}$ has a semi-regular property if there exists $\gamma_0 \> 0$ such that
\begin{align}
\displaystyle
\frac{H_{T}}{h_{T}} \leq \gamma_0 \quad \forall \mathbb{T}_h \in \{ \mathbb{T}_h \}, \quad \forall T \in \mathbb{T}_h. \label{NewGeo}
\end{align}
\end{assume}

%\subsection{Anisotropic interpolation error estimates} \label{ani=int=err}

\subsection{$L^2$-orthogonal projection}
For $T \in \mathbb{T}_h$, the $L^2$-orthogonal projection $\Pi_{{T}}^0:L^1({T}) \to \mathbb{P}^0({T})$ is defined as follows:
\begin{align*}
\displaystyle
\Pi_{{T}}^0 {\varphi} := \frac{1}{ {|{T}|_d} } \int_{{T}} {\varphi} d{x}.
\end{align*}
The following theorem provides an anisotropic error estimate for the projection $\Pi_{T}^0$.

\begin{thr} \label{thr1}
Let $p \in [1,\infty)$ and $q \in [1,\infty)$ be such that
\begin{align}
\displaystyle
W^{1,p}({T}) \hookrightarrow L^q({T}), \label{Sobolev511}
\end{align}
that is $1 - \frac{d}{p} \geq - \frac{d}{q}$. Then, for any $\hat{\varphi} \in W^{1,p}(\widehat{T})$ with ${\varphi} := \hat{\varphi} \circ {\Phi}^{-1}$,
\begin{align}
\displaystyle
\| \Pi_{T}^0 \varphi - \varphi \|_{L^q(T)} \leq c {|{T}|_d}^{\frac{1}{q} - \frac{1}{p}} \sum_{i=1}^d h_i \left\| \frac{\partial \varphi}{\partial r_i} \right\|_{L^{p}(T)}. \label{L2ortho}
\end{align}
\end{thr}

\begin{pf*}
This proof is provided in \cite[Theorem 1]{Ish23a}.
\qed	
\end{pf*}

We also define the global interpolation $\Pi_h^0$ to space $P_{dc,h}^{0}$ as
\begin{align*}
\displaystyle
(\Pi_h^0 \varphi)|_{T} := \Pi_{T}^0 (\varphi|_{T}) \quad \forall T \in \mathbb{T}_h, \quad \forall \varphi \in L^1(\Omega).
\end{align*}

\subsection{CR finite element interpolation operator}
Let ${F}_i$, $1 \leq i \leq d+1$, be the $(d-1)$-dimensional subsimplex of ${T}$ opposite to ${p}_i$. Using the barycentric coordinates $ \{ {\lambda}_i \}_{i=1}^{d+1}: \mathbb{R}^d \to \mathbb{R}$ on the reference element, the nodal basis functions are defined as follows:
\begin{align}
\displaystyle
{\theta}_i({x}) := d \left( \frac{1}{d} - {\lambda}_i({x}) \right) \quad \forall i \in \{ 1, \ldots, d+1 \}. \label{CR912}
\end{align}
The local CR interpolation operator is defined as
\begin{align}
\displaystyle
I_{{T}}^{CR}:  W^{1,1}({T})  \ni {\varphi}  \mapsto I_{{T}}^{CR} {\varphi} := \sum_{i=1}^{d+1} \left( \frac{1}{ {| {F}_i|_{d-1}} } \int_{{F}_i} {\varphi} d {s} \right) {\theta}_i \in \mathbb{P}^1({T}). \label{CR913b}
\end{align}

We then present the anisotropic CR interpolation error estimate. 

\begin{thr} \label{thr=CR}
Let $p \in [1,\infty)$ and $q \in [1,\infty)$ be such that \eqref{Sobolev511} holds true. For any $\hat{\varphi} \in W^{2,p}(\widehat{T})$ with ${\varphi} := \hat{\varphi} \circ {\Phi}^{-1}$, it holds that
\begin{align}
\displaystyle
| I_{T}^{CR} \varphi - \varphi |_{{W^{1,q}({T}) }} &\leq  c {|{T}|_d}^{ \frac{1}{q} - \frac{1}{p} }  \sum_{i=1}^d  h_i \left | \frac{\partial\varphi}{\partial r_i} \right |_{W^{1,p}(T)}. \label{CR21} 
\end{align}
\end{thr}

\begin{pf*}
This proof is found in \cite[Theorem 2]{Ish23a}.
\qed	
\end{pf*}

The vector-valued local interpolation operator 
\begin{align*}
\displaystyle
\mathcal{I}_{T}^{CR}: W^{1,1}(T)^d \to \mathbb{P}^1(T)^d
\end{align*}
is defined component-wise, that is,
\begin{align*}
\displaystyle
\mathcal{I}_{T}^{CR} v := ( I_{T}^{CR} v_1,\ldots,  I_{T}^{CR} v_d)^{\top} \quad \forall v = (v_1,\ldots,v_d)^{\top} \in W^{1,1}(T)^d.
\end{align*}
We define a global  interpolation operator $\mathcal{I}_{h}^{CR}: W^{1,1}(\Omega)^d \to V_{h}$ as
\begin{align}
\displaystyle
(\mathcal{I}_{h}^{CR} v )|_{T} = \mathcal{I}_{T}^{CR} (v |_{T})  \quad \forall T \in \mathbb{T}_h, \quad \forall v \in W^{1,1}(\Omega)^d,  \label{CR8}
\end{align}
where the space $V_h$ is defined in Section \ref{FEspaces0}.

\subsection{RT finite element interpolation operator} \label{RTsp}
For a simplex $T \subset \mathbb{R}^d$, we define the local RT polynomial space as follows:
\begin{align}
\displaystyle
\mathbb{RT}^0(T) := \mathbb{P}^0(T )^d + x \mathbb{P}^0(T ), \quad x \in \mathbb{R}^d. \label{RTp}
\end{align}
Let $\mathcal{I}_{T}^{RT^0}: W^{1,1}(T)^d \to \mathbb{RT}^0(T)$ be the RT interpolation operator such that for any $v \in W^{1,1}(T)^d$,
\begin{align}
\displaystyle
\mathcal{I}_{T}^{RT^0}: W^{1,1}(T)^d \ni v \mapsto \mathcal{I}_{T}^{RT^0} v := \sum_{i=1}^{d+1} \left(  \int_{{F}_{i}} {v} \cdot n_{F_i} d{s} \right) \theta_{T,i}^{RT^0} \in \mathbb{RT}^0(T), \label{RT4}
\end{align}
where $\theta_{T,i}^{RT^0}$ is the local shape function (e.g. \cite{ErnGue21a}) and $n_{F_i}$ is a fixed unit normal to ${F}_{i}$. 

The following two theorems are based on element $T$ satisfying  Type \roman{sone} or Type \roman{stwo} in Section \ref{two=step} when $d=3$.

\begin{thr} \label{thr3}
Let $p \in [1,\infty)$ and $T$ be the element under Conditions \ref{cond1} or \ref{cond2} and satisfy (Type \roman{sone}) in Section \ref{two=step} when $d=3$. For any $\hat{v} \in W^{1,p}(\widehat{T})^d$ with ${v} = ({v}_1,\ldots,{v}_d)^{\top} := {\Psi} \hat{v}$,
\begin{align}
\displaystyle
\| \mathcal{I}_{T}^{RT^0} v - v \|_{L^p(T)^d} 
&\leq  c \left( \frac{H_{T}}{h_{T}} \sum_{i=1}^d h_i \left \|  \frac{\partial v}{\partial r_i} \right \|_{L^p(T)^d} +  h_{T} \| \div {v} \|_{L^{p}({T})} \right). \label{RT5}
\end{align}
\end{thr}

\begin{pf*}
The proof can be found in \cite[Theorem 2]{Ish21}.
\qed
\end{pf*}

\begin{thr} \label{thr4}
Let $p \in [1,\infty)$, $d=3$, and $T$ be the element with Condition \ref{cond2} and satisfy (Type \roman{stwo}) in Section \ref{two=step}. For $\hat{v} \in W^{1,p}(\widehat{T})^3$ with ${v} = ({v}_1,v_2,{v}_3)^{\top} := {\Psi} \hat{v}$,
\begin{align}
\displaystyle
&\| \mathcal{I}_{T}^{RT^0} v - v \|_{L^p(T)^3} 
\leq c \frac{H_{T}}{h_{T}} \Biggl(  h_{T} |v|_{W^{1,p}(T)^3} \Biggr). \label{RT6}
\end{align}
\end{thr}

\begin{pf*}
The proof can be found in \cite[Theorem 3]{Ish21}.
\qed
\end{pf*}

The RT finite element space is defined as follows:
\begin{align*}
\displaystyle
V^{RT^0}_{h} &:= \{ v_h \in L^1(\Omega)^d: \  v_h|_T \in \mathbb{RT}^0(T), \ \forall T \in \mathbb{T}_h, \  [\![ v_h \cdot n ]\!]_F = 0, \ \forall F \in \mathcal{F}_h \}.
\end{align*}
We define the following global RT interpolation $\mathcal{I}_{h}^{RT^0} : V_h  \cup W^{1,1}(\Omega)^d \to V^{RT^0}_{h}$ as
\begin{align*}
\displaystyle
(\mathcal{I}_{h}^{RT^0} v )|_{T} = \mathcal{I}_{T}^{RT^0} (v|_{T}) \quad \forall T \in \mathbb{T}_h, \quad \forall v \in  V_h \cup W^{1,1}(\Omega)^d.
\end{align*}

%\color{red}
\subsection{RT interpolation errors on anisotropic elements} \label{ani=mesh=RT}
Let $T \subset \mathbb{R}^2$ be a simplex. For any $v \in H^1(T)^2$, the RT interpolation error in \eqref{RT5} with $p=2$ can be expressed as follows:
\begin{align}
\displaystyle
\| \mathcal{I}_{T}^{RT^0} v - v \|_{L^2(T)^2} 
&\leq  c \left( \frac{H_{T}}{h_{T}} \sum_{i=1}^2 h_i \left \|  \frac{\partial v}{\partial r_i} \right \|_{L^2(T)^2} +  h_{T} \| \div {v} \|_{L^{2}({T})} \right). \label{RT5=p=2}
\end{align}

The shape-regularity condition is known: there exists a constant $\gamma_1 \> 0$ such that
\begin{align*}
\displaystyle
\rho_{T} \geq \gamma_1 h_{T} \quad \forall \mathbb{T}_h \in \{ \mathbb{T}_h \}, \quad \forall T \in \mathbb{T}_h,
\end{align*}
where $\rho_T$ is the diameter of the largest ball that can be inscribed in $T$ and is equivalent to the following condition: there exists a constant $\gamma_2 \> 0$ such that for any $\mathbb{T}_h \in \{ \mathbb{T}_h \}$ and simplex $T \in \mathbb{T}_h$, we have
\begin{align*}
\displaystyle
|T|_2 \geq \gamma_2 h_{T}^2. 
\end{align*}
This proof is provided in \cite[Theorem 1]{BraKorKri08}. 

We considered the following four anisotropic elements: Let $0 \< s, \delta \ll 1$, $s,\delta \in \mathbb{R}$, and $\varepsilon \> 1$, $\varepsilon \in \mathbb{R}$.

\begin{ex}
Let $T \subset \mathbb{R}^2$ be the simplex with vertices $p_1 := (0,0)^{\top}$, $p_2 := (2s,0)^{\top}$, and $p_3 := (s , \delta s)^{\top}$. Then, we have that $h_1 = h_2 =  s \sqrt{1+\delta^2 } \approx s$, $h_T = 2s$, $|T|_2 = \delta s^2$, and
\begin{align*}
\displaystyle
\frac{h_T^2}{|T|_2} = \frac{4}{\delta} \< + \infty, \quad \frac{H_T}{h_T} = \frac{1+\delta^2}{\delta} \< + \infty.
\end{align*}
Therefore, both the shape regularity and semi-regularity conditions are satisfied. Estimate \eqref{RT5=p=2} is as follows:
\begin{align*}
\displaystyle
\| \mathcal{I}_{T}^{RT^0} v - v \|_{L^2(T)^2} 
&\leq c \left( \frac{1+\delta^2}{\delta} \sum_{i=1}^2 s \sqrt{1+\delta^2 } \left \|  \frac{\partial v}{\partial r_i} \right \|_{L^2(T)^2} +  2 s \| \div {v} \|_{L^{2}({T})} \right) \\
&\leq  c s |v|_{H^1(T)^2}.
\end{align*}
In this case, the anisotropic error \eqref{RT5=p=2} is reduced to the standard error estimate.
\end{ex}

\begin{ex}
Let $T \subset \mathbb{R}^2$ be the simplex with vertices $p_1 := (0,0)^{\top}$, $p_2 := (2s,0)^{\top}$, and $p_3 := (s ,s^{\varepsilon})^{\top}$. Then, we have that $h_1 = h_2 = \sqrt{s^2+s^{2 \varepsilon}} \approx s$, $h_T = 2 s$, $|T|_2 = s^{1+\varepsilon}$, and
\begin{align*}
\displaystyle
\frac{h_T^2}{|T|_2} = 4 s^{1-\varepsilon} \to \infty \ \text{as $s \to 0$}, \quad \frac{H_T}{h_T} = \frac{s^2 + s^{2 \varepsilon}}{s^{1+\varepsilon}} \to  \infty \ \text{as $s  \to 0$}.
\end{align*}
Therefore, neither the shape-regularity nor the semi-regularity condition is satisfied. The estimate \eqref{RT5=p=2} is as follows:
\begin{align*}
\displaystyle
\| \mathcal{I}_{T}^{RT^0} v - v \|_{L^2(T)^2} 
&\leq  c \left(  \frac{s^2 + s^{2 \varepsilon}}{s^{1+\varepsilon}} \sum_{i=1}^2 \sqrt{s^2+s^{2 \varepsilon}} \left \|  \frac{\partial v}{\partial r_i} \right \|_{L^2(T)^2} +  2s \| \div {v} \|_{L^{2}({T})} \right) \\
&\leq c s^{2-\varepsilon} |v|_{H^1(T)^2}.
\end{align*}
In this case, when $\varepsilon \> 2$, the estimate diverges as $s \to 0$.
\end{ex}

\begin{ex}
Let $T \subset \mathbb{R}^2$ be the simplex with vertices $p_1 := (0,0)^{\top}$, $p_2 := (s,0)^{\top}$, and $p_3 := (0 , \delta s)^{\top}$. Then, we have that $h_1 = s$, $h_2 = \delta s$, $h_T = s \sqrt{1 + \delta^2} \approx s$, and $|T|_2 = \frac{1}{2} \delta s^2$, and
\begin{align*}
\displaystyle
\frac{h_T^2}{|T|_2} =\frac{2(1+\delta^2)}{\delta} \< + \infty, \quad \frac{H_T}{h_T} = 2.
\end{align*}
Therefore, both the shape regularity and semi-regularity conditions are satisfied. The estimate \eqref{RT5=p=2} is as follows:
\begin{align*}
\displaystyle
\| \mathcal{I}_{T}^{RT^0} v - v \|_{L^2(T)^2} 
&\leq  c \left(  s \left \|  \frac{\partial v}{\partial r_1} \right \|_{L^2(T)^2} + \delta s \left \|  \frac{\partial v}{\partial r_2} \right \|_{L^2(T)^2}\right) \\
&\quad + c s \sqrt{1 + \delta^2} \| \div {v} \|_{L^{2}({T})}.
\end{align*}
In this case, the anisotropic error \eqref{RT5=p=2} is valid.
\end{ex}

\begin{ex}
Let $T \subset \mathbb{R}^2$ be the simplex with vertices $p_1 := (0,0)^{\top}$, $p_2 := (s,0)^{\top}$, and $p_3 := (0 ,s^{\varepsilon})^{\top}$. Subsequently, we obtain $h_1 = s$, $h_2 = s^{\varepsilon}$, $h_T = \sqrt{s^2 + s^{2 \varepsilon}} \approx s$, and $|T|_2 = \frac{1}{2} s^{1+\varepsilon}$, and
\begin{align*}
\displaystyle
\frac{h_T^2}{|T|_2} = \frac{2 (s^2 + s^{2 \varepsilon})}{s^{1+\varepsilon}} \to \infty \ \text{as $s \to 0$}, \quad \frac{H_T}{h_T} = 2.
\end{align*}
Therefore, although the shape-regularity condition is not satisfied, the semi-regularity condition is satisfied. The estimate \eqref{RT5=p=2} is as follows:
\begin{align*}
\displaystyle
\| \mathcal{I}_{T}^{RT^0} v - v \|_{L^2(T)^2} 
&\leq  c \left(  s \left \|  \frac{\partial v}{\partial r_1} \right \|_{L^2(T)^2} + s^{\varepsilon} \left \|  \frac{\partial v}{\partial r_2} \right \|_{L^2(T)^2}\right) \\
&\quad + c \sqrt{s^2 + s^{2 \varepsilon}}  \| \div {v} \|_{L^{2}({T})}.
\end{align*}
In this case, the anisotropic error \eqref{RT5=p=2} is valid.
\end{ex}

\color{black}

\subsection{Notable properties for analysis}
The following relationship holds between the RT interpolation $\mathcal{I}_{h}^{RT^0}$ and $L^2$-projection $\Pi_h^0$:

\begin{lem} \label{lem2}
It holds that
\begin{align}
\displaystyle
\div (\mathcal{I}_{T}^{RT^0} v) = \Pi_{T}^0 (\div v) \quad \forall v \in H^{1}(T)^d. \label{RT53}
\end{align}
By combining \eqref{RT53}, for any $v \in H^1(\Omega)^d$, it holds that
\begin{align}
\displaystyle
\div (\mathcal{I}_{h}^{RT^0} v) = \Pi_h^0 (\div v). \label{RT54}
\end{align}
\end{lem}

\begin{pf*}
This proof is provided in \cite[Lemma 16.2]{ErnGue21a}.
\qed
\end{pf*}

Lemma 2 presents an important relationship.

\begin{lem} \label{rel=lem2}
For any $v_h \in V^{RT^0}_{h}$ and $\psi_h \in H_0^1(\Omega) +  V_{h0}^{CR}$,
\begin{align}
\displaystyle
\sum_{T \in \mathbb{T}_h} \int_{T} ( v_h \cdot \nabla ) \psi_h dx + \sum_{T \in \mathbb{T}_h} \int_{T} \div v_h \psi_h dx  = 0. \label{RTCRre2}
\end{align}
\end{lem}

\begin{pf*}
For any $v_h \in V_h^{RT^0}$ and $\psi_h \in H_0^1(\Omega) +  V_{h0}^{CR}$, using Green’s formula and from $v_h \cdot n_F \in \mathbb{P}^0(F)$ for any $F \in \mathcal{F}_h$, we derive 
\begin{align*}
\displaystyle
&\sum_{T \in \mathbb{T}_h} \int_{T} ( v_h \cdot \nabla ) \psi_h dx + \sum_{T \in \mathbb{T}_h} \int_{T} \div v_h \psi_h dx \\
&\quad = \sum_{T \in \mathbb{T}_h} \int_{\partial T} (v_h \cdot n_{T}) \psi_h ds  \notag 
 = \sum_{F \in \mathcal{F}_h} \int_{F} [\![ (v_h  \psi_h ) \cdot n_F ]\!] ds  \notag \\
&\quad = \sum_{F \in \mathcal{F}_h^i} \int_{F} \left(  [\![ v_h \cdot n]\!]_F \{ \! \{ \psi_h \} \!\}_{\overline{\omega},F} + \{ \! \{ v_h  \} \!\}_{\omega,F} \cdot n_F [\![ \psi_h ]\!]_F \right) ds \notag \\
&\quad \quad + \sum_{F \in \mathcal{F}_h^{\partial}} \int_{F} (v_h \cdot n_F) \psi_h ds  =  0,
\end{align*}
which leads to \eqref{RTCRre2}.
\qed
\end{pf*}

\begin{note}
For any $w_h \in V_h  \cup W^{1,1}(\Omega)^d$, we use the following relation:
\begin{align*}
\displaystyle
&\sum_{T \in \mathbb{T}_h} \int_{T} ( \mathcal{I}_{h}^{RT^0} w_h \cdot \nabla ) \varphi_h dx + \sum_{T \in \mathbb{T}_h} \int_{T} (\Pi_{h}^{0} \div w_h ) \varphi_h dx  = 0, \quad \varphi_h \in V_{h0}^{CR}.
\end{align*}
\end{note}

\begin{note}
Let $T \in \mathbb{T}_h$ and $k \in \mathbb{N}$. We define a discrete space as follows:
\begin{align*}
\displaystyle
\mathbb{P}_k^{d,k} := \{q \in \mathbb{P}^k(T)^d: \  q \cdot n_T \in \mathbb{R}^k(\partial T) \},
\end{align*}
where
\begin{align*}
\displaystyle
\mathbb{R}^k(\partial T) := \{ \varphi_h \in L^2(\partial T): \ \varphi_h |_F \in \mathbb{P}^k(F) \ \forall F \in \mathcal{F}_h \}.
\end{align*}
$H(\div;\Omega)$-conforming Brezzi--Douglas--Marini finite element space (\cite{BofBreFor13}) is defined as
\begin{align*}
\displaystyle
V^{H(\div;\Omega)}_{h} &:= \{ v_h \in L^1(\Omega)^d: \  v_h|_T \in \mathbb{P}_k^{d,k}, \ \forall T \in \mathbb{T}_h, \  [\![ v_h \cdot n ]\!]_F = 0, \ \forall F \in \mathcal{F}_h \}.
\end{align*}
Let $\mathcal{I}^{H(\div;\Omega)}_h: V_h \cup W^{1,1}(\Omega)^d \to V^{H(\div;\Omega)}_{h}$ be a global interpolation operator satisfying
\begin{align*}
\displaystyle
\int_{F}  \{ \! \{ \mathcal{I}^{H(\div;\Omega)}_h v   \} \!\}_{\omega,F} \cdot n_F q_h ds = 0 \quad \forall v \in V_h \cup W^{1,1}(\Omega)^d, \quad \forall q_h \in \mathbb{R}^k(\partial T).
\end{align*}
By setting $v_h :=  \mathcal{I}^{H(\div;\Omega)}_h v$, the relation \eqref{RTCRre2} holds. Furthermore, let $V_h^c$ be an $H_0^1(\Omega)$-conforming finite element space. Let $\mathcal{L}_h: V_h \cup W^{1,1}(\Omega)^d \to V^{H(\div;\Omega)}_{h}$ be the approximate lifting operators that preserve the optimal convergence rate of the estimated error. Because $ [\![ \psi_h ]\!]_F = 0$ for any $\psi_h \in V_h^c$ and setting $v_h :=  \mathcal{L}_h v$, the relation \eqref{RTCRre2} holds.
\end{note}

\section{Discrete Sobolev inequality}
We impose a weak elliptic regularity assumption to obtain a discrete Sobolev inequality.

\begin{assume} \label{assume2}
Let  $p^{\prime} \in (\frac{2d}{d+2} , 2]$. We assume that, for any $q \in L^{p^{\prime}}(\Omega)$, the variational problem 
\begin{align}
\displaystyle
\int_{\Omega} \nabla z \cdot \nabla w dx = \int_{\Omega} q w dx \quad \forall w \in H_0^1(\Omega) \label{weakreg}
\end{align}
has a unique solution $z \in H_0^1(\Omega)$ that belongs to $W^{2,p^{\prime}}(\Omega)$. In addition, the elliptic regularity estimate
\begin{align*}
\displaystyle
\| z \|_{W^{2,p^{\prime}}(\Omega)} \leq C^{LE} \| q \|_{L^{p^{\prime}}(\Omega)}
\end{align*}
holds, where $C^{LE}$ is a positive constant that is independent of $z$.
\end{assume}

\begin{note}
From the Lax--Milgram theorem, it follows that the problem \eqref{weakreg} has a unique solution $z \in H_0^1(\Omega)$. 
\end{note}

%\color{red}
\begin{note}
In two dimensions, \cite{Dau88} may be helpful in showing that the solution $z$ belongs to $W^{2,p^{\prime}}(\Omega)$. A discontinuous Galerkin analysis of low-regularity solutions was reported by \cite{PieErn12,WihRiv11}.
\end{note}

\begin{note}
Here, we introduce an example from \cite{WihRiv11}. Let $\Omega := (0, \frac{1}{4})^2$. The function $q$ of the Poisson equation
\begin{align*}
\displaystyle
- \varDelta z  = q \quad \text{in $\Omega$}, \quad z = 0 \quad \text{on $\partial \Omega$}
\end{align*}
is provided to satisfy the exact solution:
\begin{align*}
\displaystyle
z(x,y) := x \left (x-\frac{1}{4} \right) y \left (y-\frac{1}{4} \right) r^{-2 + \alpha},
\end{align*}
where $\alpha \in (0,1]$ is a constant, and $r := (x^2 + y^2)^{\frac{1}{2}}$. Then,
\begin{align*}
\displaystyle
z \in H^1_0(\Omega) \cap W^{2,p^{\prime}}(\Omega),
\end{align*}
for all $p^{\prime} \in (1 , \frac{2}{2 - \alpha}) \subseteq (1,2)$.
\end{note}

\color{black}

\begin{lem} \label{disSov}
Suppose that Assumptions \ref{neogeo=assume} and \ref{assume2} hold. Let $p \in [2 , \infty )$ if $d=2$, or $p \in [2,6]$ if $d=3$.  Assume that there exists a positive constant $C^{dS0}$ independent of $T$ and $h$ such that
\begin{align}
\displaystyle
\max_{T \in \mathbb{T}_h} \left( {|{T}|_d}^{\frac{1}{p} - \frac{1}{2}} h_T \right) \leq C^{dS(p)}. \label{Sobmesh}
\end{align}
It then holds that
\begin{align}
\displaystyle
\| \varphi_h \|_{L^p(\Omega)} \leq C^{dS0(p)} |\varphi_h|_{H^1(\mathbb{T}_h)} \quad \forall \varphi_h \in V_{h0}^{CR}, \label{DSineq}
\end{align}
where $C^{dS0(p)}$ is a positive constant independent of $h$ and $\varphi_h$. 
\end{lem}

\begin{pf*}
Let $p^{\prime}$ be the conjugate of $p$ with $\displaystyle \frac{1}{p} + \frac{1}{p^{\prime}} = 1$; that is, $p^{\prime} \in (\frac{2d}{d+2} , 2]$. It holds that $L^{p}(\Omega) \subset L^2(\Omega) \subset L^{p^{\prime}}(\Omega)$ and for any $\xi \in L^p(\Omega)$ and $\chi \in L^2(\Omega)$,
\begin{align}
\displaystyle
\| \xi \|_{L^p(\Omega)} = \sup_{\chi \neq 0} \frac{\langle \xi , \chi \rangle_{L^p(\Omega),L^{p^{\prime}}(\Omega)}}{\| \chi \|_{L^{p^{\prime}}(\Omega)}} =  \sup_{\chi \neq 0} \frac{\int_{\Omega} \xi \chi dx}{\| \chi \|_{L^{p^{\prime}}(\Omega)}}. \label{DS=17}
\end{align}
Let $q \in L^2(\Omega)$ and $\varphi_h \in  V_{h0}^{CR}$. Consider the following problem: Find $z \in H^1_0(\Omega)$ such that
\begin{align*}
\displaystyle
- \varDelta z = q \quad \text{in $\Omega$}, \quad z=0 \quad \text{on $\partial \Omega$}.
\end{align*}
From Assumption \ref{assume2}, there exists a unique solution $z \in H^1_0(\Omega)$, the solution is in $z \in W^{2,p^{\prime}}(\Omega)$, and 
\begin{align}
\displaystyle
\| z \|_{W^{2,p^{\prime}}(\Omega)} \leq C^{LE} \| q \|_{L^{p^{\prime}}(\Omega)}. \label{DS=elliptic}
\end{align}
From the Sobolev embedding theorem, the inclusion $H^1(\Omega) \hookrightarrow L^{p}(\Omega)$ holds true. Therefore, it holds that
\begin{align*}
\displaystyle
|z|^2_{H^1(\Omega)} = \int_{\Omega} - \varDelta z z dx = \int_{\Omega} q z dx \leq \| q \|_{L^{p^{\prime}}(\Omega)} \| z \|_{L^p(\Omega)} \leq c \| q \|_{L^{p^{\prime}}(\Omega)} |z|_{H^1(\Omega)},
\end{align*}
implying
\begin{align}
\displaystyle
|z|_{H^1(\Omega)} \leq c \| q \|_{L^{p^{\prime}}(\Omega)}. \label{DS=19}
\end{align}
From \eqref{DS=17}, we have that
\begin{align}
\displaystyle
\| \varphi_h \|_{L^p(\Omega)}
&= \sup_{q \neq 0} \frac{\int_{\Omega} \varphi_h q dx}{\| q \|_{L^{p^{\prime}}(\Omega)}} = \sup_{q \neq 0} \frac{- \int_{\Omega} \varphi_h \varDelta z dx}{\| q \|_{L^{p^{\prime}}(\Omega)}} = \sup_{q \neq 0} \frac{- \int_{\Omega} \varphi_h \div(\nabla z) dx}{\| q \|_{L^{p^{\prime}}(\Omega)}}. \label{DS=20}
\end{align}
To obtain the target estimate, we first introduce the following simple calculations:
\begin{align*}
\displaystyle
- \int_{\Omega} \div (\nabla z)  \varphi_h dx
&=  \int_{\Omega} ( \Pi_h^0  \div (\nabla z) - \div (\nabla z) ) \varphi_h dx - \int_{\Omega} ( \Pi_h^0  \div (\nabla z) ) \varphi_h dx\\
&=  \int_{\Omega} ( \Pi_h^0  \div (\nabla z) - \div (\nabla z) ) ( \varphi_h - \Pi_h^0 \varphi_h) dx \\
&\quad  - \int_{\Omega} ( \div \mathcal{I}_h^{RT^0}  (\nabla z) ) \varphi_h dx\\
&= - \int_{\Omega} \div (\nabla z) \left( \varphi_h - \Pi_h^0 \varphi_h \right) dx \\
&\quad - \int_{\Omega} (\nabla z - \mathcal{I}_h^{RT^0} (\nabla z) ) \cdot \nabla_h \varphi_h dx + \int_{\Omega} \nabla z \cdot \nabla_h \varphi_h dx,
\end{align*}
where
\begin{align*}
\displaystyle
 \int_{\Omega} ( \div \mathcal{I}_h^{RT^0}  (\nabla z) ) \varphi_h dx
 &= \sum_{T \in \mathbb{T}_h} \int_{\partial T} n_T \cdot  \mathcal{I}_h^{RT^0}  (\nabla z) \varphi_h ds - \int_{\Omega}  \mathcal{I}_h^{RT^0}  (\nabla z) \cdot \nabla_h \varphi_h dx \\
 &=  \int_{\Omega} ( \nabla z -  \mathcal{I}_h^{RT^0}  (\nabla z) ) \cdot \nabla_h \varphi_h dx  - \int_{\Omega}  \nabla z \cdot \nabla_h \varphi_h dx.
\end{align*}
From the H\"older inequality, we have
\begin{align}
\displaystyle
\left| - \int_{\Omega} \div (\nabla z)  \varphi_h dx \right|
&\leq \| \varDelta z \|_{L^{p^{\prime}}(\Omega)} \|  \varphi_h - \Pi_h^0 \varphi_h \|_{L^p(\Omega)} \notag \\
&\quad + \sum_{T \in \mathbb{T}_h} \| \nabla z - \mathcal{I}_T^{RT^0} (\nabla z)  \|_{L^{p^{\prime}}(T)} \| \nabla_h \varphi_h \|_{L^p(T)} \notag \\
&\quad + |z|_{H^1(\Omega)} |\varphi_h|_{H^1(\mathbb{T}_h)} \notag \\
&=: I_1 + I_2 + I_3. \label{DS=21}
\end{align}
Because $W^{1,p}(T) \subset L^2(T)$, and using \eqref{L2ortho}, $I_1$ is estimated as
\begin{align}
\displaystyle
I_1 \leq  c \max_{T \in \mathbb{T}_h} \left( {|{T}|_d}^{\frac{1}{p} - \frac{1}{2}} h_T \right) \| q \|_{L^{p^{\prime}}(\Omega)} |\varphi_h|_{H^1(\mathbb{T}_h)}. \label{DS=22}
\end{align}
Using  the scaling argument (\cite{IshKobTsu23a}), we have
\begin{align}
\displaystyle
| \varphi_h |_{W^{1,p}(T)}
&\leq c  {|{T}|_d}^{\frac{1}{p} - \frac{1}{2}} \frac{H_T}{h_T} |\varphi_h|_{H^1(T)}. \label{DS=23}
\end{align}
For the proof, we assume that all the norms in $\mathbb{P}^1(\widehat{T})$ are equivalent. Using the semi-regular mesh condition (Assumption \ref{neogeo=assume}), \eqref{RT5}, \eqref{RT6}, \eqref{DS=elliptic}, and \eqref{DS=23}, $I_2$ is estimated as
\begin{align}
\displaystyle
I_2
&\leq c \sum_{T \in \mathbb{T}_h} {|{T}|_d}^{\frac{1}{p} - \frac{1}{2}} h_T |z|_{W^{2,p^{\prime}}(T)} |\varphi_h|_{H^1(T)} \notag \\
&\leq c \max_{T \in \mathbb{T}_h} \left( {|{T}|_d}^{\frac{1}{p} - \frac{1}{2}} h_T \right) \| q \|_{L^{p^{\prime}}(\Omega)} |\varphi_h|_{H^1(\mathbb{T}_h)}. \label{DS=25}
\end{align}
Using \eqref{DS=19}, $I_3$ is estimated as follows:
\begin{align}
\displaystyle
I_3
\leq c \| q \|_{L^{p^{\prime}}(\Omega)}  |\varphi_h|_{H^1(\mathbb{T}_h)}. \label{DS=26}
\end{align}
Using \eqref{Sobmesh}, \eqref{DS=20}, \eqref{DS=21}, \eqref{DS=22}, \eqref{DS=25}, and \eqref{DS=26}, we obtain the target inequality \eqref{DSineq}.
\qed
\end{pf*}

\begin{rem} \label{NewRem4}
Let ${|{T}|_d} \approx h_T^{d + \varepsilon}$ for any $1 \leq \varepsilon \in \mathbb{R}$. We then have
\begin{align*}
\displaystyle
 {|{T}|_d}^{\frac{1}{p} - \frac{1}{2}} h_T \approx h_T^{(d+\varepsilon)\frac{2-p}{2p} + 1}.
\end{align*}
For example, when $p=4$, if $d+\varepsilon \< 4$, there exists a positive constant $\delta \in \mathbb{R}$ such that $h_T^{\delta} \to 0$ as $h_T \to 0$. In an isotropic element, $\varepsilon = 0$, that is,  $ {|{T}|_d}^{- \frac{1}{4}} h_T \approx h_T^{-\frac{d}{4} + 1}$.
\end{rem}

%\color{red}
\begin{rem}
The proof of Lemma \ref{disSov} is based on the regularity assumption (Assumption \ref{assume2}). The regularity assumption is difficult, and the domain may be restricted. For example, when $p^{\prime} = 2$, $\Omega$ is convex. To avoid the restrictive conditions, several methods have been proposed (\cite{Bre03,PieErn12,GirLiRiv16,LasSul03}). However, these methods have been proven under the shape-regularity mesh condition. As described in Section \ref{ani=mesh=RT},  anisotropic meshes can be classified into two types: 
\begin{description}
  \item[(a)] those that include elements with large aspect ratios and those that do not satisfy the shape-regularity condition but satisfy the semi-regular condition (Assumption \ref{neogeo=assume});
  \item[(b)] those that include elements with large aspect ratios and whose partitions satisfy the shape-regularity condition.  
\end{description}
In the case of (b), it may be possible to extend these new methods; however, in the case of (a), it is difficult to apply these methods to anisotropic meshes. In \cite{Bre03}, discrete Poincar\'e inequalities for piecewise $H^1$ functions are proposed. However, the inverse, trace inequalities, and local quasi-uniformity for meshes under the shape-regular condition are used for the proof. Therefore, a careful consideration of the results in \cite{Bre03} may be necessary to eliminate the assumption that $\Omega$ is convex. For example, it may be necessary to assume the minimum angle condition for simplices within the macro elements. Further investigation of this issue remains a topic for future work.
\end{rem}

\color{black}

Using the discrete Sobolev inequality, we obtain the following stability estimates of the RT interpolation.

\begin{lem}[Stability of the RT interpolation] \label{RT=lem3}
We impose the same assumptions as in Lemma \ref{disSov}. It holds that
\begin{align}
\displaystyle
&\| v_h - \mathcal{I}_h^{RT^0} v_h \|_{L^p(\Omega)^d} \leq C^{dS1(p)} \ |v_h|_{V_h} \quad \forall v_h \in V_{h}, \label{DSRT=1} \\
&\| \mathcal{I}_h^{RT^0} v_h \|_{L^p(\Omega)^d} \leq C^{dS2(p)}  |v_h|_{V_h} \quad \forall v_h \in V_{h}, \label{DSRT=2}
\end{align}
where $C^{dS1(p)}$ and $C^{dS2(p)}$ are positive constants independent of $h$ and $\varphi_h$ but dependent on $C^{dS0(p)}$. 
\end{lem}

\begin{pf*}
For $T \in \mathbb{T}_h$, by using \eqref{RT5}, \eqref{RT6}, \eqref{Sobmesh} and \eqref{DS=23}, we have
\begin{align*}
\displaystyle
\| v_h - \mathcal{I}_T^{RT^0} v_h \|_{L^p(T)^d}
&\leq c h_T |v_h|_{W^{1,p}(T)^d}
\leq c  {|{T}|_d}^{\frac{1}{p} - \frac{1}{2}} h_T |v_h|_{H^1(T)^d},
\end{align*}
implying
\begin{align*}
\displaystyle
\| v_h - \mathcal{I}_h^{RT^0} v_h \|_{L^p(\Omega)^d} &\leq c C^{dS0(p)}  |v_h|_{V_h}.
\end{align*}
Inequality \eqref{DSRT=2} follows from the triangle inequality, \eqref{Sobmesh}, \eqref{DSineq}, and \eqref{DSRT=1}.
\qed
\end{pf*}

\section{Stability and error estimates}
\subsection{Well-posedness and the stability}
The boundness of $a_{\rot,h}$ follows from the stability of the discrete Sobolev inequality.

\begin{lem} \label{tri=lem4}
We impose the same assumptions as those in Lemma \ref{disSov} with $p=4$. It holds that for any $u(h) \in V + V_h$ and $v_h,w_h \in V_h$,
\begin{align}
\displaystyle
|a_{\rot,h}(u(h),v_h,w_h)| &\leq N^{CR} |u(h)|_{V_h}  |v_h|_{V_h}  |w_h|_{V_h}, \label{cr=3}
\end{align}
where $N^{CR}$ is a positive constant independent of $h$ but dependent on $C^{dS0(4)}$. 
\end{lem}

\begin{pf*}
From \eqref{DSRT=2} with $p=4$ and the H\"older inequality with exponents $(\frac{1}{4},\frac{1}{2},\frac{1}{4})$, the target inequality holds.
\qed
\end{pf*}

Next, we prove the well-posedness and stability of discrete problems \eqref{cr=1} and \eqref{cr=2}.

\begin{thr}
We impose the same assumptions as those in Lemma \ref{disSov}. Then, the modified CR finite element method \eqref{cr=2} has at least one solution with the stability estimate:
\begin{align}
\displaystyle
| u_h |_{V_h} \leq \frac{C^{dS2(2)}}{\nu} \| P_{\mathcal{H}} f \|_{L^2(\Omega)^d}. \label{cr=5a}
\end{align}
If the condition
\begin{align}
\displaystyle
\frac{2 N^{CR}C^{dS2(2)}}{\nu^2} \| P_{\mathcal{H}} f \|_{L^2(\Omega)^d} \< 1 \label{cr=5}
\end{align}
holds true, then problems \eqref{cr=2} and \eqref{cr=1} have unique solutions $u_h \in V_{h,\div}$ and $(u_h,p_h) \in V_{h} \times Q_h$, respectively.
\end{thr}

\begin{pf*}
For each $w_h \in V_{h,\div}$, we associate an element $R_h(w_h) \in V_{h,\div}$ as a solution to the following linear problem:
\begin{align}
\displaystyle
\nu a_{h}(R_h(w_h),v_h) + a_{\rot,h}( w_h, R_h(w_h), v_h ) = \int_{\Omega} f \cdot \mathcal{I}_{h}^{RT^0} v_h dx \quad \forall v_h \in V_{h,\div}. \label{cr=6}
\end{align}
The solvability follows from the fact that the corresponding homogeneous problem:
\begin{align}
\displaystyle
\nu a_{h}(R_h(w_h),v_h) + a_{\rot,h}( w_h, R_h(w_h), v_h ) = 0 \quad \forall v_h \in V_{h,\div} \label{cr=7}
\end{align}
has a trivial solution $R_h(w_h) = 0$. Setting $v_h := R_h(w_h)$ in \eqref{cr=7} gives
\begin{align*}
\displaystyle
0 &= \nu |R_h(w_h)|^2_{V_h} + a_{\rot,h}( w_h, R_h(w_h), R_h(w_h) ) = \nu |R_h(w_h)|^2_{V_h}.
\end{align*}
This defines the mapping $R_h: V_{h,\div} \to V_{h,\div}$. Setting $v_h := R_h(w_h)$ in \eqref{cr=6} gives
\begin{align*}
\displaystyle
| R_h(w_h) |_{V_h} \leq \frac{C^{dS2(2)}}{\nu} \| P_{\mathcal{H}} f \|_{L^2(\Omega)^d} =: R.
\end{align*}
Mapping $R_h$ maps the ball $V_{h,\div} \cap B_R := \{ v_h \in V_{h,\div}: \ |v_h|_{V_h} \leq R \}$ to itself. For any $w_h, \varphi_h \in V_{h,\div}$, the following relation is valid:
\begin{align*}
\displaystyle
0
&= \nu a_{h}(R_h(w_h) - R_h(\varphi_h),v_h) + a_{\rot,h}( w_h, R_h(w_h), v_h )  - a_{\rot,h}( \varphi_h, R_h(\varphi_h), v_h ) \\
&= \nu a_{h}(R_h(w_h) - R_h(\varphi_h),v_h) + a_{\rot,h}( w_h, R_h(w_h) - R_h(\varphi_h), v_h ) \\
&\quad + a_{\rot,h}(w_h - \varphi_h, R_h(\varphi_h), v_h ) \quad \forall v_h \in V_{h,\div}.
\end{align*}
We set $v_h := R_h(w_h) - R_h(\varphi_h)$. Using \eqref{cr=3}, we obtain
\begin{align*}
\displaystyle
\nu |R_h(w_h) - R_h(\varphi_h)|_{V_h}^2
&= - a_{\rot,h}(w_h - \varphi_h, R_h(\varphi_h), R_h(w_h) - R_h(\varphi_h) ) \\
&\leq N^{CR} |w_h - \varphi_h|_{V_h} |R_h(\varphi_h)|_{V_h}  |R_h(w_h) - R_h(\varphi_h)|_{V_h},
\end{align*}
which leads to
\begin{align*}
\displaystyle
 |R_h(w_h) - R_h(\varphi_h)|_{V_h}
 &\leq \frac{N^{CR}}{\nu} R  |w_h - \varphi_h|_{V_h}.
\end{align*}
According to Brouwer's fixed-point theorem, there exists at least one fixed point $w_h \in V_{h,\div}$ such that $w_h = R_h(w_h)$. By definition, the fixed point satisfies
\begin{align*}
\displaystyle
\nu a_{h}(w_h,v_h) + a_{\rot,h}( w_h, w_h, v_h ) = \int_{\Omega} f \cdot \mathcal{I}_{h}^{RT^0} v_h dx \quad \forall v_h \in V_{h,\div},
\end{align*}
and 
\begin{align*}
\displaystyle
| w_h |_{V_h} \leq \frac{C^{dS2(2)}}{\nu} \| P_{\mathcal{H}} f \|_{L^2(\Omega)^d}.
\end{align*}
Let $u_h^1,u_h^2 \in V_{h,\div}$ be solutions to \eqref{cr=2}. We set $w_h := u_h^1 - u_h^2$ . Using \eqref{cr=4}, we obtain
\begin{align*}
\displaystyle
\nu a_{h}(w_h,v_h) 
&= a_{\rot,h}( u_h^2, u_h^2, v_h )  - a_{\rot,h}( u_h^1, u_h^1, v_h ) \\
&= - a_{\rot,h}( w_h, u_h^2, v_h ) - a_{\rot,h}( u_h^1, w_h, v_h ).
\end{align*}
Setting $v_h := w_h$ yields
\begin{align*}
\displaystyle
\nu |w_h|_{V_h}^2
&\leq N^{CR} \left( |u_h^2|_{V_h} + |u_h^1|_{V_h} \right )  |w_h|_{V_h}^2.
\end{align*}
If $\nu \> 2  N^{CR} R$, $ |w_h|_{V_h}^2 \leq 0$ holds, leading to $w_h \equiv 0$.

For the solution $u_h \in V_{h,\div}$ of \eqref{cr=2}, we show that there exists a corresponding pressure $p _h \in Q_h$ such that
\begin{align}
\displaystyle
b_h(v_h , p_h)
&= \int_{\Omega} f \cdot \mathcal{I}_h^{RT^0} v_h dx - \nu a_h(u_h,v_h) - a_{\rot,h}(u_h,  u_h,  v_h ) \quad \forall v_h \in V_{h}. \label{cr=8}
\end{align}
The right-hand side of \eqref{cr=8} represents a linear function $\ell(\cdot)$ on $V_h$ with $\ell(v_h) = 0$, $v_h \in V_{h,\div}$. The problem \eqref{cr=1} is well-posed from the discrete inf-sup conditions \eqref{cr=4} and \eqref{cr=8}.
\qed
\end{pf*}

\subsection{Error estimates}
We introduce a Jensen-type inequality \cite[Exercise 12.1]{ErnGue21a}. Let $r,s$ be two nonnegative real numbers and $\{ x_i \}_{i \in I}$ be a finite sequence of nonnegative numbers. It then holds that
\begin{align}
\displaystyle
\begin{cases}
\left( \sum_{i \in I} x_i^s \right)^{\frac{1}{s}} \leq \left( \sum_{i \in I} x_i^r \right)^{\frac{1}{r}} \quad \text{if $r \leq s$},\\
\left( \sum_{i \in I} x_i^s \right)^{\frac{1}{s}} \leq \card(I)^{\frac{r-s}{rs}} \left( \sum_{i \in I} x_i^r \right)^{\frac{1}{r}} \quad \text{if $r \> s$}.
\end{cases} \label{jensen}
\end{align}

For the analysis, we define trilinear forms $a_{\rot,h}^*: V \times V \times V_h \to \mathbb{R}$ and $\varepsilon_{\rot,h}: V \times V \times V_h \to \mathbb{R}$ as
\begin{align}
\displaystyle
a_{\rot,h}^*(u, v, w_h ) &:= \sum_{T \in \mathbb{T}_h} \int_T (\curl u) \times v \cdot \mathcal{I}_h^{RT^0} w_h dx \label{rot*} \\
&= \sum_{T \in \mathbb{T}_h} \int_T \left\{  (v \cdot \nabla) u \cdot  \mathcal{I}_h^{RT^0} w_h - (\mathcal{I}_h^{RT^0} w_h \cdot \nabla) u \cdot  v  \right\}dx, \notag \\
\varepsilon_{\rot,h}(u, v, w_h ) &:= a_{\rot,h}(u, v, w_h ) - a_{\rot,h}^*(u, v, w_h ). \notag
\end{align}

The boundness of $a_{\rot,h}^*$ and $\varepsilon_{\rot,h}$ follows from the stability of the discrete Sobolev inequality.

\begin{lem} \label{tri=lem5}
We impose the same assumptions as those in Lemma \ref{disSov} with $p=4$. It holds that for any $u,v \in V$ and $w_h \in V_{h}$,
\begin{align}
\displaystyle
|a_{\rot,h}^*(u,v, w_h)| 
&\leq N^{CR*} |u|_{V}  |v |_{V}  | w_h |_{V_h}, \label{cr=3*} \\
|\varepsilon_{\rot,h}(u, v ,w_h)| 
&\leq N^{CR**} |u |_{V}  \|  \mathcal{I}_h^{RT^0} v - v \|_{L^4(\Omega)^d}  | w_h |_{V_h}, \label{cr=3***} 
\end{align}
where $N^{CR*}$ and $N^{CR**}$ are positive constants independent of $h$ but dependent on $C^{dS0(4)}$.
\end{lem}

\begin{pf*}
From \eqref{DSRT=2} with $p=4$ and the H\"older inequality with exponents $(\frac{1}{4},\frac{1}{2},\frac{1}{4})$, the inequality \eqref{cr=3*} holds.

Given
\begin{align*}
\displaystyle
&\varepsilon_{\rot,h}(u, v, w_h ) \\
&= \sum_{T \in \mathbb{T}_h} \int_T \left \{ ( (\mathcal{I}_h^{RT^0} v - v ) \cdot \nabla) u \cdot \mathcal{I}_h^{RT^0} w_h - (\mathcal{I}_h^{RT^0} w_h \cdot \nabla )u \cdot (\mathcal{I}_h^{RT^0} v - v) \right \} dx,
\end{align*}
using \eqref{DSRT=2} with $p=4$, the H\"older inequality with exponents $(\frac{1}{4},\frac{1}{2},\frac{1}{4})$ yields the inequality \eqref{cr=3***}.
\qed
\end{pf*}

We first introduce two lemmata for error analysis.

\begin{lem} \label{lin=lem6}
 Let $T$ be the element that satisfies Condition \ref{cond1} or \ref{cond2} and is of Type \roman{sone}, as described in Section \ref{two=step}, when $d=3$. We impose the same assumptions as those in Lemma \ref{disSov} with $p=4$. Let $(u,p) \in V \times Q$ be the solution to \eqref{ns1} and $u_h \in V_{h,\div}$ be the solution to \eqref{cr=2}. We assume $u \in H^2(\Omega)$. For any $\varphi_h \in V_{h,\div}$, it holds that
\begin{align}
\displaystyle
&| a_{\rot,h}(u_h, u_h, \varphi_h ) - a_{\rot,h}^*(u,u, \varphi_h)| \notag\\
&\quad \leq N^{CR} |u_h - u + \varphi_h|_{V_h} |u_h|_{V_h} |\varphi_h|_{V_h}  + N^{CR} |\varphi_h|_{V_h}^2 |u_h|_{V_h} \notag\\
&\quad \quad + N^{CR}|u|_V |u_h - u + \varphi_h|_{V_h} |\varphi_h|_{V_h} \notag\\
&\quad \quad + N_1^{CR**} |u |_{V} \sum_{j=1}^d {\sum_{T \in \mathbb{T}_h}} h_j \left \|  \frac{\partial u}{\partial r_j} \right \|_{L^4(T)^d}   | \varphi_h |_{V_h},  \label{cr=10}
\end{align}
where $ N_1^{CR**}$ is a positive constant independent of $h$ but dependent on $ N^{CR**}$.
\end{lem}

\begin{pf*}
For any $\varphi_h \in V_{h,\div}$, as $a_{\rot,h}(u, \varphi_h, \varphi_h) = 0$, it holds that
\begin{align*}
\displaystyle
&a_{\rot,h}(u_h, u_h, \varphi_h ) - a_{\rot,h}^*(u,u, \varphi_h) \\
&\quad = a_{\rot,h}(u_h - u, u_h, \varphi_h ) + a_{\rot,h}(u, u_h, \varphi_h )  -  a_{\rot,h}(u, u, \varphi_h )  \\
&\quad \quad +  a_{\rot,h}(u, u, \varphi_h ) - a_{\rot,h}^*(u,u, \varphi_h) \\
&\quad = a_{\rot,h}(u_h - u, u_h, \varphi_h ) + a_{\rot,h}(u, u_h - u, \varphi_h ) + \varepsilon_{\rot,h}(u,u, \varphi_h)  \\
&\quad = a_{\rot,h}(u_h - u + \varphi_h, u_h, \varphi_h ) -  a_{\rot,h}( \varphi_h, u_h, \varphi_h ) \\
&\quad \quad  + a_{\rot,h}(u, u_h - u + \varphi_h, \varphi_h ) + \varepsilon_{\rot,h}(u,u, \varphi_h).
\end{align*}
Using \eqref{RT5}, \eqref{cr=3}, \eqref{cr=3*}, \eqref{cr=3***}, and $\div u = 0$, we obtain
\begin{align*}
\displaystyle
&|a_{\rot,h}(u_h, u_h, \varphi_h ) - a_{\rot,h}^*(u,u, \varphi_h) | \\
&\quad \leq N^{CR} |u_h - u + \varphi_h|_{V_h} |u_h|_{V_h} |\varphi_h|_{V_h}  + N^{CR} |\varphi_h|_{V_h}^2 |u_h|_{V_h}\\
&\quad \quad + N^{CR}|u|_V |u_h - u + \varphi_h|_{V_h} |\varphi_h|_{V_h} \\
&\quad \quad + N^{CR**} |u |_{V}  \|  \mathcal{I}_h^{RT^0} u - u \|_{L^4(\Omega)^d}  | \varphi_h |_{V_h} \\
&\quad \leq N^{CR} |u_h - u + \varphi_h|_{V_h} |u_h|_{V_h} |\varphi_h|_{V_h}  + N^{CR} |\varphi_h|_{V_h}^2 |u_h|_{V_h}\\
&\quad \quad + N^{CR}|u|_V |u_h - u + \varphi_h|_{V_h} |\varphi_h|_{V_h} \\
&\quad \quad + N_1^{CR**} |u |_{V} \sum_{j=1}^d {\sum_{T \in \mathbb{T}_h} }h_j \left \|  \frac{\partial u}{\partial r_j} \right \|_{L^4(T)^d}   | \varphi_h |_{V_h},
\end{align*}
which is the target inequality.
\qed	
\end{pf*}

\begin{lem} \label{lem7}
Let $T$ be the element that satisfies Condition \ref{cond1} or \ref{cond2} and is of Type \roman{sone}, as described in Section \ref{two=step}, when $d=3$. We now impose Assumption \ref{neogeo=assume}. Let $(u,p) \in V \times Q$ and $u_h \in V_{h,\div}$  be the solutions to \eqref{ns1} and \eqref{cr=2}, respectively. We assume $(u,p) \in H^2(\Omega) \times H^1(\Omega)$. For any $\varphi_h \in V_{h,\div}$, it holds that
\begin{align}
\displaystyle
&\left| \nu a_h(u,\varphi_h)  + a_{\rot,h}^*(u,u, \varphi_h)  - \int_{\Omega} f \cdot \mathcal{I}_{h}^{RT^0} \varphi_h dx \right| \notag\\
&\quad \leq  c \nu \left( \sum_{j=1}^d  \sum_{T \in \mathbb{T}_h} h_j \left |  \frac{\partial u}{\partial r_j} \right |_{H^1(T)^d} + h \| \varDelta u \|_{L^2(\Omega)^d} \right) |\varphi_h|_{V_h}. \label{cr=15}
\end{align}
\end{lem}

\begin{pf*}
Let $\varphi_h \in V_{h,\div}$. From $f = - \nu \varDelta u + (\curl u) \times u + \nabla p $, \eqref{rot*}, and $\int_{\Omega} \nabla p \cdot \mathcal{I}_{h}^{RT^0} \varphi_h dx = 0$, we obtain
\begin{align*}
\displaystyle
a_{\rot,h}^*(u,u,  \varphi_h)  - \int_{\Omega} f \cdot \mathcal{I}_{h}^{RT^0} \varphi_h dx = \nu \int_{\Omega} \varDelta u \cdot \mathcal{I}_{h}^{RT^0} \varphi_h dx.
\end{align*}
Using the relation between the RT and CR finite elements (see Lemma \ref{rel=lem2}) and \eqref{RT54} yields
\begin{align*}
\displaystyle
\nu a_h(u,\varphi_h)  &= \nu \sum_{i=1}^d \sum_{T \in \mathbb{T}_h} \int_{T}( \nabla u_{i} - \mathcal{I}_{h}^{RT^0} \nabla u_i) \cdot \nabla \varphi_{h,i} dx - \nu \int_{\Omega} \Pi_h^0 \varDelta u \cdot \varphi_h dx.
\end{align*}
Therefore, from the stability of $\Pi_h^0$, \eqref{L2ortho}, and \eqref{RT5}, it holds that 
\begin{align*}
\displaystyle
&\left| \nu a_h(u,\varphi_h)  + a_{\rot,h}^*(u,u, \varphi_h)  - \int_{\Omega} f \cdot \mathcal{I}_{h}^{RT^0} \varphi_h dx \right| \notag\\
&\ \leq \nu \sum_{i=1}^d \sum_{T \in \mathbb{T}_h} \int_{T} | \nabla u_{i} - \mathcal{I}_{h}^{RT^0} \nabla u_i|   | \nabla \varphi_{h,i}| dx + \nu \int_{\Omega} |\varDelta u| |\mathcal{I}_{h}^{RT^0} \varphi_h - \varphi_h| dx  \notag \\
&\ + \nu \int_{\Omega} | \varDelta u - \Pi_h^0 \varDelta u |  | \varphi_h - \Pi_h^0 \varphi_h | dx  \notag \\
&\leq c \nu \left( \sum_{i,j=1}^d  \sum_{T \in \mathbb{T}_h} h_j \left \|  \frac{\partial \nabla u_i}{\partial r_j} \right \|_{L^2(T)^d} + h \| \varDelta u \|_{L^2(\Omega)^d} \right) |\varphi_h|_{V_h},
\end{align*}
which leads to the target estimate using the Jensen-type inequality \eqref{jensen}.
\qed
\end{pf*}

The following theorem is an error estimate for the velocity:

\begin{thr} \label{main=th7}
 Let $T$ be the element satisfying Condition \ref{cond1} or \ref{cond2} and be of Type \roman{sone}, as discussed in Section \ref{two=step}, when $d=3$. We impose the same assumptions as those in Lemma \ref{disSov} with $p=4$. Let $(u,p) \in V \times Q$ and $u_h \in V_{h,\div}$ be the solutions of \eqref{ns1} and \eqref{cr=2}, respectively. We assume $(u,p) \in H^2(\Omega)^d \times H^1(\Omega)$. Furthermore, we impose uniqueness conditions in \eqref{ns6} and \eqref{cr=5}. It then holds that
\begin{align}
\displaystyle
|u-u_h|_{V_h}
&\leq C(\sigma,\| P_{\mathcal{H}} f \|_{L^2(\Omega)^d})  \notag\\
&\hspace{-0.8cm} \times \Biggl(  \sum_{j=1}^d  \sum_{T \in \mathbb{T}_h} h_j \left |  \frac{\partial u}{\partial r_j} \right |_{H^1(T)^d} + h \| \varDelta u \|_{L^2(\Omega)^d}   + \sum_{j=1}^d { \sum_{T \in \mathbb{T}_h} } h_j \left \|  \frac{\partial u}{\partial r_j} \right \|_{L^4(T)^d} \Biggr), \label{cr=16}
\end{align}
where $C(\sigma,\| P_{\mathcal{H}} f \|_{L^2(\Omega)^d})$ is a positive constant independent of $h$ and $\sigma$ is defined in the proof.
\end{thr}

\begin{pf*}
For any $w_h \in V_{h,\div}$, we have that
\begin{align}
\displaystyle
|u-u_h|_{V_h} &\leq  |u-w_h|_{V_h}  +  |w_h-u_h|_{V_h}. \label{cr=17}
\end{align}
Using \eqref{ns7}, \eqref{cr=5a}, \eqref{cr=10} with $\varphi_h := w_h-u_h$, and \eqref{cr=15} yields
\begin{align*}
\displaystyle
 \nu |\varphi_h|_{V_h}^2
&= \nu a_h(\varphi_h , \varphi_h) \\
&\hspace{-0.8cm}  = \nu a_h(w_h - u, \varphi_h) + \nu a_h(u,\varphi_h)  + a_{\rot,h}^*(u,u, \varphi_h)  - \int_{\Omega} f \cdot \mathcal{I}_{h}^{RT^0} \varphi_h dx \\
&\hspace{-0.8cm}  \quad + a_{\rot,h}(u_h, u_h, \varphi_h ) - a_{\rot,h}^*(u,u,   \varphi_h) \\
&\hspace{-0.8cm}  \leq \nu |u-w_h|_{V_h} |\varphi_h|_{V_h} + \left| \nu a_h(u,\varphi_h)  + a_{\rot,h}^*(u,u,  \varphi_h)  - \int_{\Omega} f \cdot \mathcal{I}_{h}^{RT^0} \varphi_h dx \right| \\
&\hspace{-0.8cm}  \quad + \left| a_{\rot,h}(u_h, u_h, \varphi_h ) - a_{\rot,h}^*(u,u, \varphi_h) \right|,
\end{align*}
which leads to
\begin{align}
\displaystyle
 |\varphi_h|_{V_h}
 &\leq |u-w_h|_{V_h} + c \left( \sum_{j=1}^d  \sum_{T \in \mathbb{T}_h} h_j \left |  \frac{\partial u}{\partial r_j} \right |_{H^1(T)^d} + h \| \varDelta u \|_{L^2(\Omega)^d} \right) \notag \\
  &\hspace{-1.0cm} + \frac{N^{CR} C^{dS2(2)}}{\nu^2} \| P_{\mathcal{H}} f \|_{L^2(\Omega)^d} |w_h - u|_{V_h} + \frac{N^{CR} C^{dS2(2)}}{\nu^2} \| P_{\mathcal{H}} f \|_{L^2(\Omega)^d} |\varphi_h|_{V_h} \notag \\
   &\hspace{-1.0cm} + \frac{N^{CR} C_P}{\nu^2} \| P_{\mathcal{H}} f \|_{L^2(\Omega)^d} |w_h - u |_{V_h} + \frac{N_1^{CR**} C_P}{\nu^2}  \| P_{\mathcal{H}} f \|_{L^2(\Omega)^d} \sum_{j=1}^d { \sum_{T \in \mathbb{T}_h}}  h_j \left \|  \frac{\partial u}{\partial r_j} \right \|_{L^4(T)^d}.
\end{align}
Condition \eqref{cr=5} leads to
\begin{align*}
\displaystyle
1 - \frac{N^{CR}C^{dS2(2)}}{\nu^2} \| P_{\mathcal{H}} f \|_{L^2(\Omega)^d}
= \frac{\nu^2 - N^{CR}C^{dS2(2)} \| P_{\mathcal{H}} f \|_{L^2(\Omega)^d}}{\nu^2} \> 0.
\end{align*}
Thus, $\sigma := \nu^2 - N^{CR}C^{dS2(2)} \| P_{\mathcal{H}} f \|_{L^2(\Omega)^d} \> 0$. Therefore, using \eqref{cr=17}, \eqref{ns6}, and \eqref{cr=5} yields
\begin{align*}
\displaystyle
|u-u_h|_{V_h} &\leq \left(2+ C_0(\sigma,\| P_{\mathcal{H}} f \|_{L^2(\Omega)^d}) \right) |u-w_h|_{V_h} \\
&\quad + C_1(\sigma,\| P_{\mathcal{H}} f \|_{L^2(\Omega)^d}) \sum_{j=1}^d { \sum_{T \in \mathbb{T}_h} } h_j \left \|  \frac{\partial u}{\partial r_j} \right \|_{L^4(T)^d} \\
&\quad + C_2(\sigma)  \left( \sum_{j=1}^d  \sum_{T \in \mathbb{T}_h} h_j \left |  \frac{\partial u}{\partial r_j} \right |_{H^1(T)^d} + h \| \varDelta u \|_{L^2(\Omega)^d} \right).
\end{align*}
As $w_h \in V_{h,\div}$ is set arbitrarily, setting $w_h := \mathcal{I}_{T}^{CR} u$ yields the target inequality with the Jensen-type inequality.
\qed
\end{pf*}

Next, we prove the $ L^2$ error estimate for pressure.

\begin{thr}
Let $T$ be the element satisfying Condition \ref{cond1} or \ref{cond2} and be of Type \roman{sone}, as discussed in Section \ref{two=step}, when $d=3$. We impose the same assumptions as those in Lemma \ref{disSov} with $p=4$. Let $(u,p) \in V \times Q$ and $(u_h,p_h) \in V_{h} \times Q_h$ be the solutions of \eqref{ns1} and \eqref{cr=1}, respectively. We assume $(u,p) \in H^2(\Omega) \times H^1(\Omega)$. Furthermore, we impose the uniqueness conditions \eqref{ns6} and \eqref{cr=5}. It then holds that
\begin{align}
\displaystyle
\|p-p_h\|_{Q_h}
&\leq C(\sigma)  \Biggl( \sum_{j=1}^d  \sum_{T \in \mathbb{T}_h} h_j \left |  \frac{\partial u}{\partial r_j} \right |_{H^1(T)^d} + h \| \varDelta u \|_{L^2(\Omega)^d}  \notag\\
&\quad +  |u |_{V} \sum_{j=1}^d { \sum_{T \in \mathbb{T}_h} } h_j \left \|  \frac{\partial u}{\partial r_j} \right \|_{L^4(T)^d} + \sum_{j=1}^d \sum_{T \in \mathbb{T}_h}  h_j \left\| \frac{\partial p}{\partial r_j} \right\|_{L^{2}(T)} \Biggr), \label{cr=19}
\end{align}
where $C(\sigma)$ is a positive constant independent of $h$.
\end{thr}

\begin{pf*}
To estimate the pressure error $ \| p - p_h \|_{Q_h}$, we use the inf-sup stability relation \eqref{cr=4}.  For an arbitrary $q_h \in Q_h$, it follows that
\begin{align}
\displaystyle
\| p - p_h \|_{Q_h}
&\leq  \| p - q_h \|_{Q_h} + \| q_h - p_h \|_{Q_h} \notag \\
&\leq  \| p - q_h \|_{Q_h} + \frac{1}{\beta_0}  \sup_{v_h \in V_{h}} \frac{b_h(v_h, q_h - p_h)}{|v_h|_{V_{h}}} \notag \\
&\leq  \| p - q_h \|_{Q_h} \notag \\
&\quad + \frac{1}{\beta_0}  \sup_{v_h \in V_{h}} \frac{b_h(v_h, q_h - p)}{|v_h|_{V_{h}}} + \frac{1}{\beta_0}  \sup_{v_h \in V_{h}} \frac{b_h(v_h, p - p_h)}{|v_h|_{V_{h}}} \notag \\
&\leq \left( 1 + \frac{1}{\beta_0} \right)  \| p - q_h \|_{Q_h} + \frac{1}{\beta_0}  \sup_{v_h \in V_{h}} \frac{b_h(v_h, p - p_h)}{|v_h|_{V_{h}}}. \label{cr=20}
\end{align}
Let $v_h \in V_{h}$. From $f = - \nu \varDelta u + (\curl u) \times u + \nabla p $ and \eqref{rot*}, we obtain that
\begin{align}
\displaystyle
&- \int_{\Omega} f \cdot \mathcal{I}_{h}^{RT^0} \varphi_h dx \notag \\
&\quad = \nu \int_{\Omega} \varDelta u \cdot \mathcal{I}_{h}^{RT^0} \varphi_h dx - a_{\rot,h}^*(u,u, \varphi_h) -\int_{\Omega} \nabla p \cdot \mathcal{I}_h^{RT^0} v_h dx. \label{cr=21}
\end{align}
Using the relation between the RT and CR finite elements (see Lemma \ref{rel=lem2}) and \eqref{RT53} yields
\begin{align}
\displaystyle
\nu a_h(u,v_h)  &= \nu \sum_{i=1}^d \sum_{T \in \mathbb{T}_h} \int_{T}( \nabla u_{i} - \mathcal{I}_{h}^{RT^0} \nabla u_i) \cdot \nabla v_{h,i} dx - \nu \int_{\Omega} \Pi_h^0 \varDelta u \cdot v_h dx. \label{cr=22}
\end{align}
Using the Gauss--Green formula and \eqref{RT54}, we obtain
\begin{align}
\displaystyle
- \int_{\Omega} \nabla p \cdot \mathcal{I}_h^{RT^0} v_h dx
&= \int_{\Omega} p \Pi_h^0 \divh v_h dx = \int_{\Omega} p \divh v_h dx. \label{cr=23}
\end{align}
Here, we used $\Pi_h^0 \divh v_h =  \divh v_h$ because $\div v_h$ is constant on $T$. Using \eqref{cr=21}, \eqref{cr=22}, and \eqref{cr=23}, we obtain
\begin{align}
\displaystyle
b_h(v_h, p - p_h)
&=  b_h(v_h, p ) - \int_{\Omega} f \cdot \mathcal{I}_h^{RT^0} v_h dx + \nu a_h(u_h,v_h) + a_{\rot,h}(u_h,  u_h,  v_h ) \notag \\
&= \nu a_h(u_h - u,v_h) + a_{\rot,h}(u_h,  u_h,  v_h ) - a_{\rot,h}^*(u, u, v_h ) \notag \\
&\quad + \nu \int_{\Omega} \varDelta u \cdot \mathcal{I}_h^{RT^0} v_h dx - \nu \int_{\Omega} \Pi_h^0 \varDelta u \cdot v_h dx  \notag\\
&\quad + \nu \sum_{i=1}^d \sum_{T \in \mathbb{T}_h} \int_{T}( \nabla u_{i} - \mathcal{I}_{h}^{RT^0} \nabla u_i) \cdot \nabla v_{h,i} dx \notag \\
&=: I_1+I_2+I_3+I_4+I_5+I_6. \label{cr=24}
\end{align}
 $I_1$ is easily estimated as follows:
\begin{align}
\displaystyle
|I_1| &\leq \nu |u - u_h|_{V_h} |v_h|_{V_h}. \label{cr=25}
\end{align}
Using \eqref{ns6}, \eqref{ns7}, \eqref{cr=5a}, \eqref{cr=5}, \eqref{RT5}, \eqref{cr=3}, \eqref{cr=3*}, \eqref{cr=3***}, and $\div u = 0$, the estimate of $I_2+I_3$ follows from
\begin{align}
\displaystyle
|I_2+I_3|
&\leq | a_{\rot,h}(u_h - u , u_h, v_h ) | 
+ | a_{\rot,h}(u, u_h - u , v_h ) | + | \varepsilon_{\rot,h}(u,u, v_h) | \notag \\
&\leq  N^{CR} |u_h - u |_{V_h} |u_h|_{V_h} |v_h|_{V_h}  + N^{CR}|u|_V |u_h - u |_{V_h} |v_h|_{V_h} \notag\\
&\quad + N_1^{CR**} |u |_{V} \sum_{j=1}^d { \sum_{T \in \mathbb{T}_h} } h_j \left \|  \frac{\partial u}{\partial r_j} \right \|_{L^4(T)^d} | v_h |_{V_h}. \label{cr=26}
\end{align}
As proof of Lemma \ref{lem7}, $I_4+I_5+I_6$ is estimated as
\begin{align}
\displaystyle
|I_4+I_5 + I_6| \leq c \nu \left( \sum_{j=1}^d  \sum_{T \in \mathbb{T}_h} h_j \left |  \frac{\partial u}{\partial r_j} \right |_{H^1(T)^d} + h \| \varDelta u \|_{L^2(\Omega)^d} \right) |v_h|_{V_h}. \label{cr=27}
\end{align}
From \eqref{L2ortho}, it holds that
\begin{align}
\displaystyle
\inf_{q_h \in Q_h} \| p - q_h \|_{Q_h} 
&\leq \| p - \Pi_{h}^0 p \|_{Q_h} \leq c \sum_{j=1}^d \sum_{T \in \mathbb{T}_h} h_j \left\| \frac{\partial p}{\partial r_j} \right\|_{L^{2}(T)}. \label{cr=28}
\end{align}
By using \eqref{cr=20}, \eqref{cr=24}, \eqref{cr=25}, \eqref{cr=26}, \eqref{cr=27}, and \eqref{cr=28} with \eqref{cr=16}, we obtain the target estimate.
\qed
\end{pf*}

\section{Numerical experiments}
This section presents the results of the numerical experiments. Let $d=2$ and $\Omega := (0,1)^2$. The function $f$ of the Navier--Stokes equation \eqref{intro1b},
\begin{align}
\displaystyle
- \nu \varDelta u + (\curl u) \times u + \nabla p = f \quad \text{in $\Omega$}, \quad \div u  = 0 \quad \text{in $\Omega$}, \quad u = g \quad \text{on $\partial \Omega$}, \label{num=1}
\end{align}
is given as satisfying exact solutions $(u,p)$, where $g:\Omega \to \mathbb{R}^2$ is a given function and $g \in H^{\frac{1}{2}}(\partial \Omega)^2$ satisfies $\int_{\partial \Omega} g \cdot n ds = 0$. We apply the Picard iteration to linearise the nonlinear problem \eqref{cr=1}. For $n=0,1,\ldots$, find $(u_h^{n+1} , p_h^{n+1}) \in V_h \times Q_h$ such that, for any $(v_h,q_h) \in V_h \times Q_h$,
\begin{subequations} \label{nume=1}
\begin{align}
\displaystyle
\nu a_h(u_h^{n+1},v_h) + a_{\rot,h}(u_h^{n},  u_h^{n+1},  v_h ) + b_h(v_h , p_h^{n+1})
&= \int_{\Omega} f \cdot \mathcal{I}_h^{RT^0} v_h dx, \label{nume=1a} \\
b_h(u_h^{n+1} , q_h) &= 0, \label{nume=1b}
\end{align}
\end{subequations}
{with the initial data $u_h^0 := \mathcal{I}_h^{CR} u_{ex}$, $p_h^0 := \Pi_h^0 p_{ex}$, where $u_{ex}$ and $p_{ex}$ are exact solutions.} The boundary condition is imposed by $u_h^{n+1}|_{\partial \Omega} = g_h$ for an appropriate approximation $g_h$ of $g$. We proceed with the iteration until the following end criterion is satisfied:
\begin{align*}
\displaystyle
|u_h^{n+1} - u_h^{n} |_{V_h} + \| p_h^{n+1} - p_h^n \|_{Q_h} \< 10^{-10} \left( | u_h^{n} |_{V_h} + \| p_h^n \|_{Q_h} \right).
\end{align*}
If an exact solution $u$ is known, then the errors ${e_N := u - u_N}$ and ${e_{2N} := u - u_{2N}}$ are computed numerically for the {two division numbers $N$ and $2N$.} The convergence indicator $r$ is defined as follows:
\begin{align*}
\displaystyle
r = \frac{1}{\log(2)} \log \left( \frac{\| e_N \|_X}{\| e_{2N} \|_X} \right).
\end{align*}
We compute the convergence order with respect to the norms defined by
\begin{align*}
\displaystyle
Err(V_h) &:= \frac{| u - u_N |_{V_h}}{|u|_{V_h}}, \quad Err(L^2):= \frac{\| u - u_N \|_{L^2(\Omega)^d}}{\|u\|_{L^2(\Omega)^d}},\\
Err(Q_h) &:= \frac{\| p - p_N \|_{Q_h}}{\|p\|_{Q_h}}.
\end{align*}
Numerical calculations are performed to confirm that the mesh conditions are satisfied. The semiregularity mesh condition defined in \eqref{NewGeo} is equivalent to the maximum angle condition. The following parameters are computed:
\begin{align*}
\displaystyle
\textit{MinAngle} := \max_{T \in \mathbb{T}_h} \frac{|L_3|^2}{{|{T}|_d}}, \quad \textit{MaxAngle} := \max_{T \in \mathbb{T}_h} \frac{|L_1| |L_2|}{{|{T}|_d}},
\end{align*}
where $L_i$, $i=1,2,3,$ denote the edges of the simplex $T \in \mathbb{T}_h$ with $|L_1| \leq |L_2| \leq |L_3|$. Furthermore, to verify the condition \eqref{Sobmesh} with $p=4$, we compute
\begin{align*}
\displaystyle
\textit{DisSov} :=  \max_{T \in \mathbb{T}_h} \left( {|{T}|_d}^{- \frac{1}{4}} h_T \right).
\end{align*}
See Tables \ref{table001}, \ref{table002} and \ref{table003}.

$\# Np$ denotes the number of nodal points on $V_{h} \times Q_h$, including those on the boundary. We used the GMRES(500) method with the ILU(0) preconditioner for the calculations.

\subsection{Example 1}
We set $\varphi (x_1,x_2) := 64 x_1^2(x_1-1)^2x_2^2(x_2-1)^2$. The function $f$ of the Navier--Stokes equation \eqref{num=1} with $\nu = 10^{-1}$ and $g = 0$ is given such that the exact solution is
\begin{align*}
\displaystyle
\begin{pmatrix}
 u_1  \\
  u_2
\end{pmatrix}
&= \curl \varphi =
\begin{pmatrix}
 \frac{\partial \varphi}{\partial x_2}  \\
  - \frac{\partial \varphi}{\partial x_1} 
\end{pmatrix},\\
p &= \frac{1}{2}(u_1^2 + u_2^2) -0.1238397581254773 + 10^5 (1 - x_2)^3 - \frac{10^5}{4}.
\end{align*}
It is noteworthy that $\int_{\Omega} p dx \approx 0$. 

\begin{figure}[htbp]
   \includegraphics[keepaspectratio, scale=0.15]{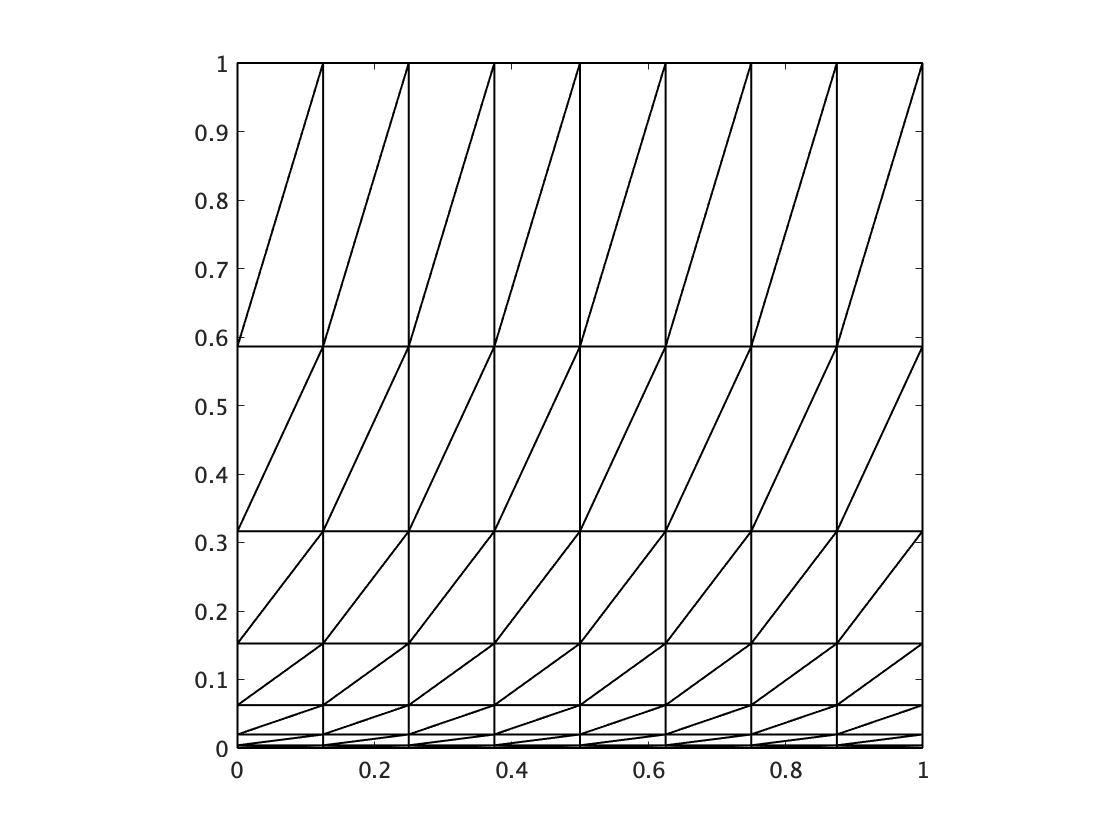}
    \caption{(Mesh \Roman{lone}) Anisotropic mesh, $N=8$ and $\varepsilon = 2.0$}
     \label{fig00}
\end{figure}

Let $N \in \{ 4,8,16,32,64,128\}$ be the division number of each side of the bottom and height edges of $\Omega$. The following mesh partitions are considered.  Let $(x_1^i, x_2^i)^T$ and $i \in \mathbb{N}$ be grip points of the triangulations $\mathbb{T}_h$ defined as follows: 
\begin{description}
  \item[(Mesh \Roman{lone})] We use the anisotropic mesh partitions.
\begin{align*}
\displaystyle
x_1^i := \frac{i}{N}, \quad x_2^i := \left ( \frac{i}{N} \right)^{\varepsilon}, \quad  i \in \{0, \ldots, N \},
\end{align*}
where $\varepsilon \in \{ 1.0, 2.0, 4.0 \}$ (Fig. \ref{fig00}). Condition \eqref{Sobmesh} is satisfied when $\varepsilon = 1.0, 2.0$ but not when $\varepsilon = 4.0$.
\end{description}

\begin{table}[h]
\caption{Mesh (\Roman{lone}), mesh conditions, $\varepsilon = 2.0$}
\centering
\begin{tabular}{l|l|l|l|l} \hline
$N$ & $\# Np$&  $\textit{MinAngle} $ & $\textit{MaxAngle}$ & $\textit{DisSov}$ \\ \hline \hline
4 & 144 & 8.50     & 2.00  & 1.04199   \\
8& 544 & 1.63e+01   & 2.00 &  7.63521e-01   \\
16& 2,112 & 3.21e+01   &  2.00 &  5.95764e-01 \\
32& 8,320 & 6.41e+01  & 2.00 & 5.00244e-01   \\
64 &  33,024 & 1.28e+02   & 2.00 & 4.20500e-01   \\
128 & 131,584 & 2.56e+02   &  2.00 & 3.53564e-01    \\
\hline
\end{tabular}
\label{table001}
\end{table}

\begin{table}[h]
\caption{Mesh (\Roman{lone}), mesh conditions, $\varepsilon = 4.0$}
\centering
\begin{tabular}{l|l|l|l|l} \hline
$N$ & $\# Np$&  $\textit{MinAngle} $ & $\textit{MaxAngle}$ & $\textit{DisSov}$ \\ \hline \hline
4 & 144 &   1.28031e+02   & 2.00  & 1.68200   \\
8& 544 &  1.02400e+03   & 2.00 & 2.00000    \\
16& 2,112 &8.19200e+03   &  2.00 &   2.37841   \\
32& 8,320 & 6.55360e+04 & 2.00 &  2.82843   \\
64 &  33,024 & 5.24288e+05  & 2.00 &   3.36359   \\
128 & 131,584 &  4.19430e+06   &  2.00 &  4.00000   \\
\hline
\end{tabular}
\label{table002}
\end{table}

The theoretical results of the modified CR finite element method for an anisotropic mesh partition are numerically verified and presented in Tables \ref{table3}, \ref{table4}, and \ref{table5}. When $\epsilon = 4.0$, the condition \eqref{Sobmesh} is not satisfied. However, the numerical result indicates an optimal convergence rate as the meshes become finer. 

\begin{table}[h]
\caption{Mesh (\Roman{lone}), Example 1, $\varepsilon = 1.0$}
\centering
\begin{tabular}{l|l|l|l|l|l|l|l} \hline
$N$ &  $h $ & $Err(V_h)$ & $r$ & $Err(L^2)$  & $r$ & $Err(Q_h)$ & $r$ \\ \hline \hline
4 & 3.54e-01        & 9.30891e-01   &  &  5.57356e-01  &   &  2.77363e-01  &  \\
8 & 1.77e-01  &    5.06405e-01  &  0.88& 1.63541e-01  & 1.77  & 1.39270e-01  & 0.99  \\
16 & 8.84e-02 &  2.59214e-01   & 0.97  &  4.33267e-02  & 1.92  &  6.97005e-02  & 1.00   \\
32 &  4.42e-02  &  1.30439e-01  &0.99&  1.10344e-02  & 1.97  & 3.48582e-02  & 1.00 \\
64 &    2.21e-02    & 6.53276e-02   &  1.00  & 2.77257e-03   & 1.99  & 1.74301e-02  & 1.00  \\
128 & 1.10e-02  &  3.26775e-02   & 1.00 & 6.93973e-04   &  2.00  &  8.71516e-03  & 1.00 \\
\hline
\end{tabular}
\label{table3}
\end{table}

\begin{table}[h]
\caption{Mesh (\Roman{lone}), Example 1, $\varepsilon = 2.0$}
\centering
\begin{tabular}{l|l|l|l|l|l|l|l} \hline
$N$ &  $h $ & $Err(V_h)$ & $r$ & $Err(L^2)$  & $r$ & $Err(Q_h)$ & $r$ \\ \hline \hline
4 & 5.04e-01       & 1.04386   &  & 7.54616e-01    &   & 2.28331e-01  &  \\
8 & 2.66e-01 &  6.00986e-01  & 0.80 &  2.50020e-01   &  1.59 &  1.13984e-01  & 1.00 \\
16 & 1.36e-01&  3.14178e-01  &0.94& 7.08474e-02    &  1.82 & 5.69444e-02   & 1.00   \\
32 &  6.90e-02 & 1.59284e-01  & 0.98  & 1.85985e-02  &1.93   &  2.84658e-02    & 1.00 \\
64 &   3.47e-02   &7.99483e-02    & 0.99 & 4.71970e-03   &   1.98 &  1.42321e-02  & 1.00  \\
128 & 1.74e-02  &  4.00138e-02   & 1.00 &   1.18479e-03   &  1.99  & 7.11597e-03   & 1.00 \\
\hline
\end{tabular}
\label{table4}
\end{table}

\begin{table}[h]
\caption{Mesh (\Roman{lone}), Example 1, $\varepsilon = 4.0$}
\centering
\begin{tabular}{l|l|l|l|l|l|l|l} \hline
$N$ &  $h $ & $Err(V_h)$ & $r$ & $Err(L^2)$  & $r$ & $Err(Q_h)$ & $r$ \\ \hline \hline
4 &  7.28e-01  & 1.13521    &  &  9.15578e-01 &   &  3.45283e-01  &  \\
8 & 4.32e-01   &  8.34160e-01   & 0.44  & 5.29158e-01  &  0.79   & 1.65246e-01  & 1.06  \\
16 &2.36e-01   &  4.72051e-01  & 0.82  & 1.80204e-01  & 1.55 & 8.17474e-02  &  1.02   \\
32 &1.23e-01   & 2.47274e-01   & 0.93 & 5.25128e-02  &  1.78 &  4.07539e-02 & 1.00 \\
64 &  6.30e-02  &  1.25537e-01  & 0.98 &1.39353e-02   & 1.91  & 2.03619e-02  &  1.00 \\
128 &3.19e-02   &  6.30344e-02  &  0.99 &  3.54646e-03  &1.97& 1.01790e-02  & 1.00 \\
\hline
\end{tabular}
\label{table5}
\end{table}

\subsection{Example 2}
The second example demonstrates the robustness of the proposed scheme against large irrotational body forces on anisotropic mesh partitions. We consider the numerical example inspired by \cite{LinMer16,QuiPie20,YanHeZha22}. The exact solution $(u,p)$ to \eqref{num=1} with velocity and pressure is given by
\begin{align*}
\displaystyle
\begin{pmatrix}
 u_1  \\
  u_2
\end{pmatrix}
&=
\begin{pmatrix}
 - \left( x_2 - \frac{1}{2} \right)  \\
x_1 - \frac{1}{2}
\end{pmatrix}, \\
p &= \left(x_1- \frac{1}{2} \right)^2  + \left(x_2- \frac{1}{2} \right)^2  - \frac{1}{6} + 10^5 (1 - x_2)^3 - \frac{10^5}{4}.
\end{align*}
It is noteworthy that $\int_{\Omega} p dx = 0$. By setting $\nu = 1$, the force is exactly irrotational, that is, 
\begin{align*}
\displaystyle
\begin{pmatrix}
 f_1  \\
 f_2
\end{pmatrix}
&=
\begin{pmatrix}
0 \\
- 3 \times 10^5 (1 - x_2)^2
\end{pmatrix}.
\end{align*}

\begin{figure}[htbp]
   \includegraphics[keepaspectratio, scale=0.15]{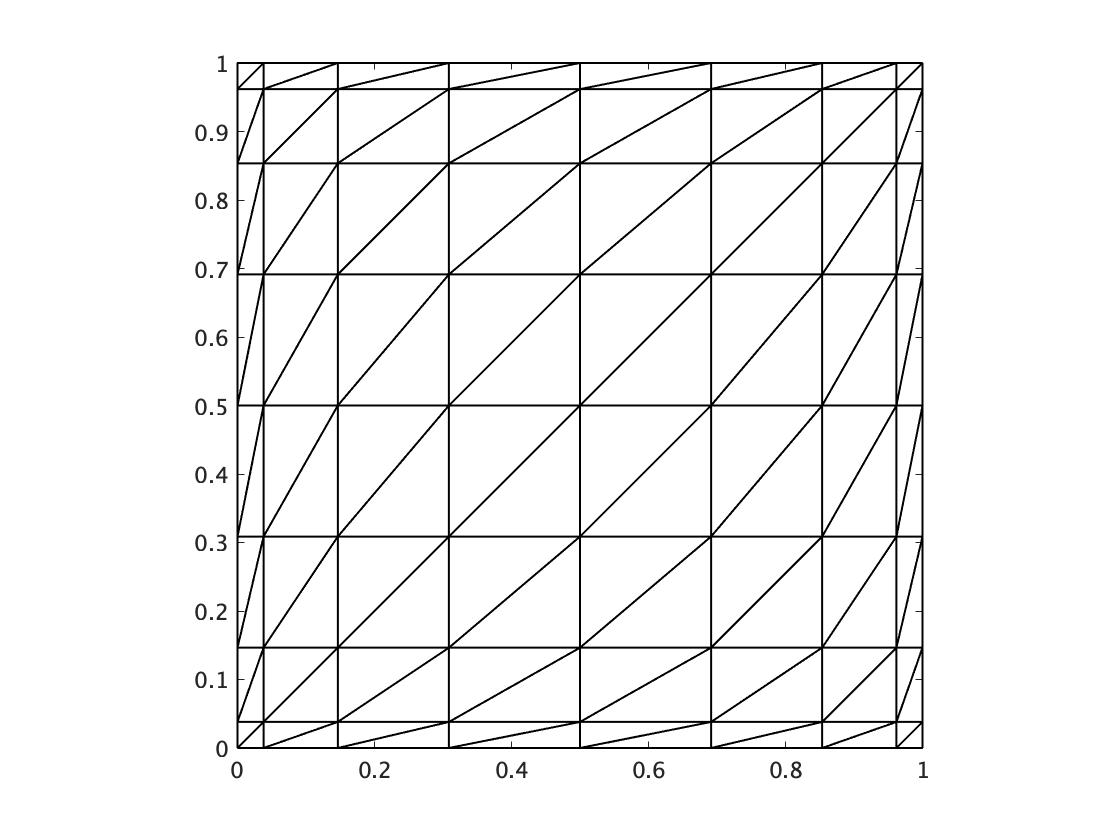}
    \caption{(Mesh \Roman{ltwo}) Anisotropic mesh, $N=8$}
     \label{fig002}
\end{figure}

Let $N \in \{ 4,8,16,32,64,128\}$ be the division number of each side of the bottom and height edges of $\Omega$. The following mesh partitions are considered. Let $(x_1^i, x_2^i)^T$ be grip points of the triangulations $\mathbb{T}_h$ defined as follows. Let $i \in \mathbb{N}$. 
\begin{description}
  \item[(Mesh \Roman{ltwo})  Anisotropic mesh which comes from \cite{CheLiuQia10} (Fig. \ref{fig002})]
\begin{align*}
\displaystyle
x_1^i := \frac{1}{2}\left( 1 - \cos \left( \frac{i \pi}{N} \right) \right), \quad x_2^i := \frac{1}{2}\left( 1 - \cos \left( \frac{i \pi}{N} \right) \right), \quad  i \in \{0, \ldots, N \}.
\end{align*}
\end{description}

\begin{table}[h]
\caption{Mesh (\Roman{ltwo}), mesh conditions}
\centering
\begin{tabular}{l|l|l|l|l} \hline
$N$ & $\# Np$&  $\textit{MinAngle} $ & $\textit{MaxAngle}$ & $\textit{DisSov}$ \\ \hline \hline
4 & 144 &  5.65685    & 2.00  &  1.00000      \\
8& 544 & 1.04525e+01    & 2.00 & 7.94187e-01    \\
16& 2,112 &  2.05033e+01  &  2.00 & 6.66204e-01   \\
32& 8,320 &  4.08092e+01   & 2.00 &  5.59870e-01   \\
64 &  33,024 &   8.15201e+01   & 2.00 &   4.70722e-01   \\
128 & 131,584 &    1.62991e+02  &  2.00 & 3.95813e-01   \\
\hline
\end{tabular}
\label{table003}
\end{table}

Tables \ref{table6}, and \ref{table7} numerically verify the theoretical results of the modified CR finite element method for an anisotropic mesh partition and confirm the robustness towards irrotational body forces on anisotropic mesh partitions. 

\begin{table}[h]
\caption{Mesh (\Roman{lone}), Example 2, $\varepsilon = 1.0$}
\centering
\begin{tabular}{l|l|l|l|l|l|l|l} \hline
$N$ &  $h $ & $Err(V_h)$ & $r$ & $Err(L^2)$  & $r$ & $Err(Q_h)$ & $r$ \\ \hline \hline
4 & 3.54e-01  & 9.09364e-07   &  & 5.47195e-07 &   &  2.77362e-01  &  \\
8 & 1.77e-01  & 2.66354e-06   &  - &1.24705e-06   & -  & 1.39270e-01   & 0.99  \\
16 & 8.84e-02 & 1.97022e-06   & - & 1.24596e-06  &  - &  6.97007e-02   &   1.00  \\
32 &  4.42e-02  & 1.73889e-06   & -  &9.04173e-07   & -   &   3.48583e-02  & 1.00 \\
64 &    2.21e-02   &   1.26862e-06   & - &  5.57509e-07& -   & 1.74301e-02   & 1.00  \\
128 & 1.10e-02  &  1.43621e-06   & - &  8.86565e-07  &  - &  8.71518e-03   &  1.00  \\
\hline
\end{tabular}
\label{table6}
\end{table}

\begin{table}[h]
\caption{Mesh (\Roman{ltwo}), Example 2}
\centering
\begin{tabular}{l|l|l|l|l|l|l|l} \hline
$N$ &  $h $ & $Err(V_h)$ & $r$ & $Err(L^2)$  & $r$ & $Err(Q_h)$ & $r$ \\ \hline \hline
4 & 5.00e-01   & 2.98226e-06   &  & 1.08150e-06  &   &  2.87956e-01  &  \\
8 & 2.71e-01   &   2.81107e-06   & -  &1.70024e-06   & -  &  1.49758e-01  &   0.94 \\
16 &  1.38e-01  &   4.52069e-06  & - &  2.75827e-06   & -  & 7.54093e-02  &  0.99   \\
32 & 6.93e-02   &  2.36901e-06  & - & 9.65821e-07 &  - &   3.77670e-02 & 1.00 \\
64 & 3.47e-02    &   2.73752e-06  & -  &  1.11624e-06  & -   & 1.88912e-02   &  1.00  \\
128 &  1.74e-02  &  2.08281e-06  &  - &  8.56957e-07  &  - &   9.44656e-03 &  1.00  \\
\hline
\end{tabular}
\label{table7}
\end{table}

\section{Concluding remarks} \label{sec=conc}
We developed a modified CR finite element method for the rotation form of the stationary incompressible Navier--Stokes equation on anisotropic mesh partitions. The proposed method achieves pressure-independent velocity error estimates. Numerical examples were used to verify the theoretical results. In the future, we will extend the proposed approach to time-dependent Navier--Stokes problems and high-Reynolds-number flows. 

Moreover, in Numerical Example 2, we imposed an inhomogeneous Dirichlet boundary condition on the simply connected domain, although we did not study this case. In general, we assume that $\Omega \subset \mathbb{R}^d$, $d \in \{ 2,3\}$, is open,  bounded, and connected, but multiply-connected, and its boundary $\Gamma$ is Lipschitz-continuous for the analysis. We consider $\Gamma_0$ as the exterior boundary of $\Omega$, and the other components of $\Gamma$ using $\Gamma_i$, $1 \leq i \leq M$. For $M \in \mathbb{N}$, we assume that
\begin{align}
\displaystyle
\int_{\Gamma_i} g \cdot n ds = 0, \quad 0 \leq i \leq M, \label{SOC}
\end{align}
where $n = (n_1,\ldots,n_d)^T$ denotes the unit outwards normal to $\Gamma_i$, $0 \leq i \leq M$. \eqref{SOC} is called the stringent outflow condition. In addition, from the divergence formula, 
\begin{align}
\displaystyle
 \int_{\Gamma} g \cdot n ds = \int_{\Omega} \div u dx = 0. \label{GOC}
\end{align}
This is known as the general outflow condition. We plan to extend the proposed approach to problems imposed by other boundary conditions, including the inhomogeneous Dirichlet boundary conditions.

\begin{acknowledgements}
We thank the Reviewers for their valuable comments and suggestions.
\end{acknowledgements}

\end{document}